\documentclass[11pt]{article}

\usepackage{amsmath,amssymb,amsthm,mathtools,amsfonts,xcolor,hyperref,enumitem, bbm}
\usepackage{mathrsfs}
\usepackage{authblk}
\usepackage{algorithm}
\usepackage{algpseudocode}
\usepackage{graphicx, float}
\usepackage{graphicx,xcolor}
\usepackage{subcaption}
\usepackage[margin=1in]{geometry}
\usepackage{setspace}
\numberwithin{equation}{section}

\newtheorem{theorem}{Theorem}[section]
\newtheorem{assumption}[theorem]{Assumption}
\newtheorem{lemma}[theorem]{Lemma}
\newtheorem{corollary}[theorem]{Corollary}
\newtheorem{definition}[theorem]{Definition}
\newtheorem{proposition}[theorem]{Proposition}
\theoremstyle{remark}
\newtheorem{remark}[theorem]{Remark}

\newcommand{\calK}{\mathcal K}

\renewcommand{\P}{\mathbb{P}}

\newcommand{\J}{\mathcal J}

\newcommand{\dd}{\mathrm d}

\newcommand{\calB}{\mathcal B}

\newcommand{\calX}{C^\infty_\diamond(\mathbb T^d)}

\newcommand{\pair}[2]{\left\langle #1,#2\right\rangle}

\DeclareMathOperator*{\argmin}{arg\,min}

\providecommand{\keywords}[1]
{
  \small	
  \textbf{\textbf{Keywords---}} #1
}

\usepackage{geometry}
\usepackage{verbatim}
\renewcommand{\P}{\mathbb{P}}

\theoremstyle{plain}

\title{Mirror flow for semimartingale optimal transport}

\author[1]{François Buet-Golfouse}
\author[1,2]{Anaïs Després}
\author[3]{Anna Kazeykina}
\author[2]{Zhenjie Ren}

\affil[1]{AIML Global Markets, Barclays}
\affil[2]{LaMME, Université Évry Paris-Saclay, Évry-Courcouronnes, France}
\affil[3]{LMO, Université Paris-Saclay, Orsay, France}
\begin{document}
\maketitle

\begin{abstract}
Semimartingale optimal transport arises in a wide range of applications, including mathematical finance and flow-based generative models. A recurring computational challenge is that the dual energy functional is often poorly conditioned with respect to the dual potential, leading to slow convergence of gradient descent algorithms. To address this issue, we investigate a mirror flow approach that effectively preconditions the dual energy functional. We establish the geometric convergence of the  mirror descent sequence. To obtain a practical numerical algorithm, we introduce a feasible discrete-time scheme based on first-order linearization. We prove that the scheme is well-defined and establish a uniform-in-time error bound for the accumulated linearization error, under a suitable smallness condition on the initial error. Finally, we apply the proposed mirror flow scheme to the specific-entropic martingale optimal transport (SEMOT) problem \cite{BuetGolfouseEtAl2026SEMOT}  and observe a substantial improvement in computational efficiency.
\end{abstract}


\section{Introduction}

Semimartingale optimal transport (SMOT) seeks the least costly
stochastic dynamics connecting a prescribed initial law $\mu$ to a
prescribed terminal law $\nu$~\cite{TaTo13}. We consider continuous
semimartingales on $[0,1]$ with dynamics
$\mathrm dX_t=b_t^P\,\mathrm dt+\sigma_t^P\,\mathrm dW_t$
under a candidate law $P$, where $b_t^P$ is the drift and
$2a_t^P=\sigma_t^P(\sigma_t^P)^\top$ is the instantaneous covariance.
Let $\mathcal P(\mu,\nu)$ be the set of admissible laws with
$X_0\sim\mu$ and $X_1\sim\nu$. The problem is
\[
V(\nu)=\inf_{P\in\mathcal P(\mu,\nu)}\mathcal C(P),
\qquad
\mathcal C(P)=\mathbb E^P\!\left[
\int_0^1 L(t,X_t,a_t^P,b_t^P)\,\mathrm dt
\right].
\]
The running cost $L$ can penalize both drift and diffusion along the
path. Applications include calibrating financial models to option
prices~\cite{GuoLoeperWang2022} and inferring cell trajectories from
population snapshots~\cite{LavenantEtAl2024}.

The terminal constraint makes the primal problem difficult to solve
directly, so we relax it with a terminal potential $\psi$, which assigns
a cost to each endpoint. For fixed $\psi$, we minimize
$\mathcal C(P)+\mathbb E^P[\psi(X_1)]$ over laws with initial law $\mu$
and no terminal constraint. Let $\mathcal A(\psi)$ be the resulting value
and $T(\psi)$ the terminal law of the optimal process. The dual objective
is
\[
D_\nu(\psi)=\mathcal A(\psi)-\int\psi\,\mathrm d\nu .
\]
It is concave and satisfies $D_\nu(\psi)\le V(\nu)$, with equality at the
supremum under suitable conditions~\cite{TaTo13}. Its gradient is
$T(\psi)-\nu$, the mismatch between the terminal law produced by $\psi$
and the target, so $\psi$ maximizes $D_\nu$ exactly when $T(\psi)=\nu$.
Each evaluation of $T(\psi)$ requires solving one stochastic control
problem.

Existing numerical methods for SMOT follow two main approaches. The first
is PDE-based and works with the dual formulation, updating the potential
by gradient or quasi-Newton methods. Finite-difference approximations and
projected gradient updates are studied in~\cite{TaTo13}, and HJB-based
gradient and quasi-Newton procedures in~\cite{GuoLangreneLoeperNing2022}.
Related implementations solve nonlinear HJB equations to calibrate
local-stochastic-volatility, joint SPX--VIX and equity--interest-rate
models~\cite{GuoLoeperWang2022,GuoLoeperOblojWang2022,
JosephLoeperObloj2024,JosephLoeperObloj2026}, and early PDE formulations
of martingale transport lead to augmented-Lagrangian and ADMM
methods~\cite{GuoLoeperWang2019}. The second approach uses discrete
transport formulations. For martingale optimal transport at a fixed
set of dates, convergent linear programming approximations are obtained in~\cite{GuoObloj2019} by discretizing the marginals and relaxing the martingale
constraint. For continuous-time SMOT, discretizing
physical time yields a multimarginal problem that can be solved by
Sinkhorn iterations~\cite{BenamouChazareixLoeper2024}. This replaces
repeated backward HJB solves with updates on the time grid, but grid
refinement increases the number of marginals.

Gradient methods for the dual face a conditioning difficulty. The terminal
law can respond much more strongly to some perturbations of $\psi$ than to
others of the same size, so the Hessian of $D_\nu$ can be ill conditioned,
and a step size small enough to keep the most sensitive directions stable
yields slow progress in the others. Quasi-Newton methods address this by
approximating the Hessian from the history of gradients and iterates. We
instead obtain the current derivative of the terminal map directly by
linearizing the HJB--Fokker--Planck system, so that, at the discrete
level, our updates use the Hessian of the minimized objective \(-D_\nu\), represented by \(-DT_\psi\). This requires
additional linear solves but gives a preconditioner adapted to the current
control problem. Our goal is to use this information to accelerate dual
optimization within the backward--forward PDE framework, without the
growth in the number of marginals caused by time refinement in
Sinkhorn-type methods.

Our approach is based on mirror flow. Given a convex functional $\Phi$,
called the \emph{mirror potential}, and an objective $F$ to be minimized,
mirror flow evolves the coordinate $D\Phi(\psi)$ along the negative
derivative of $F$. In optimization time $s$,
\[
\frac{\mathrm d}{\mathrm ds}D\Phi(\psi_s)=-DF(\psi_s),
\qquad
\dot\psi_s=-[D^2\Phi(\psi_s)]^{-1}DF(\psi_s),
\]
whenever $D^2\Phi(\psi_s)$ is invertible. The inverse Hessian rescales the
gradient according to the local geometry generated by $\Phi$, so mirror
flow is the continuous-time counterpart of preconditioned gradient
descent. When $D^2\Phi$ approximates $D^2F$, the flow approximates a
Newton flow.

For SMOT we minimize $F=-D_\nu$ and take the negative penalized control
value $\Phi=-\mathcal A$ as mirror potential. Then $D\Phi(\psi)=-T(\psi)$
and, under the assumptions of Section~\ref{sec:main}, the Hessian is
$D^2\Phi(\psi)=-DT_\psi$. Since $D^2F=D^2\Phi$, the flow
$\dot\psi_s=(-DT_{\psi_s})^{-1}\bigl(T(\psi_s)-\nu\bigr)$
is the continuous Newton flow for the dual. In the mirror coordinate
$p_s=T(\psi_s)$ it becomes linear,
\[
\frac{\mathrm dp_s}{\mathrm ds}=\nu-p_s,
\qquad
p_s-\nu=e^{-s}(p_0-\nu),
\]
so the terminal mismatch decays exponentially wherever the flow is
defined. The underlying geometry is the Bregman geometry of $\Phi$, with
divergence
$\calB_\Phi(\psi'\mid\psi)
=\Phi(\psi')-\Phi(\psi)-\langle D\Phi(\psi),\psi'-\psi\rangle\ge 0$.
It is the natural geometry in which to study the dual problem and its
convergence: if $\psi^\star$ satisfies $T(\psi^\star)=\nu$, then
$D_\nu(\psi^\star)-D_\nu(\psi)=\calB_\Phi(\psi\mid\psi^\star)$, so the dual gap
is exactly the Bregman distance from $\psi$ to the solution. Mirror flow is
steepest descent in this geometry, whereas the Euclidean metric underlying
ordinary gradient methods ignores how the terminal law responds to the
potential. The main difficulty is to show that the infinite-dimensional
Newton-type step is well defined and that its linearization can be
controlled.

Mirror and second-order methods have been studied in several transport
settings, but for different variables and dynamics. In classical optimal
transport, Sinkhorn has been interpreted as mirror descent and mirror
flow~\cite{Leger2021,AubinFrankowskiKorbaLeger2022,%
KarimiHsiehKrause2024}, while truncated Newton methods approximate
curvature corrections for a discrete transport dual%
~\cite{KemertasFarahmandJepson2025}. In stochastic transport, exponential
convergence has been proved for a gradient flow of the Bass functional in
a one-dimensional martingale setting%
~\cite{BackhoffPammerSchachermayer2025}, and affine mirror dynamics have
been studied on path laws with controlled drift and fixed
diffusion~\cite{BaiChaintronRenWang2026}. Exponential convergence has also
been established for Sinkhorn iterations in entropic martingale optimal
transport (EMOT)~\cite{ChenConfortiRenWang2024,KazeykinaRenWang2026} and
for SISTA, a method combining Sinkhorn updates with gradient
steps~\cite{BifronteLeySylvestre2026}. These results concern alternating
updates for static martingale couplings, under assumptions on compactly
supported marginals and the cost. Our scheme instead updates the terminal
potential of a continuous-time stochastic control problem.

Our differentiability analysis relies on the stability of HJB solutions
with respect to terminal data. Differentiability with respect to finite-dimensional families of terminal data and local strong convexity of the associated value map are studied in~\cite{ChassagneuxCrepeyLoeperRuelland2026}. Smooth dependence of the potentials on the marginal data for the multimarginal Schr\"odinger system is established in~\cite{CarlierLaborde2020}. Our analysis additionally concerns the terminal
marginal map of the coupled HJB--Fokker--Planck system and, in particular,
the invertibility of its derivative in Sobolev and H\"older spaces.

Our main contributions are:
\begin{enumerate}
\item \emph{Mirror flow for SMOT.} We introduce a mirror flow on terminal
  potentials for SMOT with controlled drift and diffusion. For the exact
  implicit discretization with step size $\tau$,
  $D\Phi(\bar\psi_{k+1})-D\Phi(\bar\psi_k)=-\tau\,DF(\bar\psi_{k+1})$, the terminal
  laws $\bar p_k=T(\bar\psi_k)$ satisfy $\bar p_{k+1}-\nu=(\bar p_k-\nu)/(1+\tau)$. We prove
  geometric decay of the terminal mismatch and of the dual objective gap,
  provided the iterates exist and some terminal potential produces the
  target law $\nu$.

\item \emph{Differentiability and invertibility of $T$.} For autonomous
  costs on the torus, we prove Fr\'echet differentiability of $T$ and
  identify its derivative through the linearized HJB--Fokker--Planck
  system, with a quadratic Taylor remainder. We first establish the weak invertibility
of $DT_\psi$ from $H^1_\diamond$ onto
$H^{-1}_\diamond$. Under additional
smoothness assumptions, pseudodifferential
elliptic regularity yields a strong
H\"older inverse, which in turn ensures
well-definedness of the linearized
recursion for smooth data.

\item \emph{Error of the linearized scheme.} We replace the nonlinear
  mirror steps by first-order linearized updates in which the desired
  changes in terminal density are prescribed in advance by the exact
  mirror iteration, so errors accumulated in earlier iterations are not
  corrected. Under a local smallness condition, we prove
that the terminal densities produced by this scheme remain within
$O(\tau)$, in $C^{m,\alpha}$, of the
prescribed mirror-density path.

\item \emph{Numerical illustration.} We apply the linearized scheme to the
  one-dimensional SEMOT benchmark of~\cite{BuetGolfouseEtAl2026SEMOT},
  where the discrete dual Hessian is poorly conditioned (condition number
  approximately $92.2$ initially and $863.2$ at the final mirror
  iterate). The scheme reaches the prescribed terminal accuracy with
  approximately 100 times fewer updates than ordinary gradient ascent on
  the same grid, a speedup of approximately 35 in a single run.
\end{enumerate}

The remainder of the paper is organized as follows.
Section~\ref{sec:Preliminaries} introduces the SMOT problem, its dual
formulation and the notation. Section~\ref{sec:main} presents the mirror
flow and our main results on convergence, differentiability and
invertibility of the terminal map, and the error of the linearized scheme.
Section~\ref{sec:numerical-experiments} applies the linearized scheme to
the SEMOT benchmark of~\cite{BuetGolfouseEtAl2026SEMOT} and compares it
with ordinary gradient ascent. Section~\ref{sec:proofs-mirror-terminal} contains the proofs of the mirror-flow
results, the differentiability of the terminal map,
and the weak invertibility of its derivative.
Section~\ref{sec:pseudodifferential-strong-inverse} develops the pseudodifferential
analysis needed to establish strong
H\"older invertibility.
Finally, Section~\ref{sec:proofs-linearized} proves well-definedness
of the linearized recursion and the
uniform-in-time tracking estimate.

\section{Preliminaries}
\label{sec:Preliminaries}

\paragraph{Notation.}
Fix $d,m\in\mathbb N$, $m\geq2$, and $\alpha\in(0,1)$.
We work on $\mathbb T^d:=\mathbb R^d/\mathbb Z^d$, equipped with
normalized Lebesgue measure, and set $Q:=[0,1]\times\mathbb T^d$.
Spatial derivatives are denoted by $\partial_i$, $\partial_{ij}$,
$\nabla$, and $\nabla^2$; repeated spatial indices are summed from
$1$ to $d$.
Let $\mathbb S^d$ be the symmetric $d\times d$ matrices,
$\mathbb S^d_+$ their positive-semidefinite cone, and
$\mathbb E^d:=\mathbb S^d\times\mathbb R^d$,
$\mathbb E^d_+:=\mathbb S^d_+\times\mathbb R^d$.
We use the Frobenius product $A:B:=\operatorname{Tr}(A^\top B)$
and the joint Euclidean structure
\[
(A,b)\boldsymbol\cdot(B,c):=A:B+b\cdot c,
\qquad
|(A,b)|^2:=A:A+|b|^2.
\]
Matrix inequalities refer to the positive-semidefinite order,
and $I=I_d$ is the identity matrix. For a scalar function $v$, set
$\J v:=(\nabla^2v,\nabla v)\in\mathbb E^d$.

Let $\mathcal D'(\mathbb T^d)$ denote the space of distributions on $\mathbb T^d$ and
$\langle\cdot,\cdot\rangle$ the relevant duality pairing,
written in either order and agreeing with integration whenever
applicable. For probability measures, we write
$P_n\Longrightarrow P$ for weak convergence,
meaning that
$
\int f\,dP_n\longrightarrow\int f\,dP
$
for every bounded continuous function $f$.

Define
\[
\mathcal Y:=\{r\in\mathcal D'(\mathbb T^d):\langle r,1\rangle=0\},
\qquad
F_\diamond:=F\cap\mathcal Y.
\]
We use this notation for $F=C^{q,\beta},C^\infty,H^1,H^{-1}$,
where $q\in\mathbb N_0:=\{0,1,\ldots\}$ and $\beta\in(0,1)$.
To work with functions that are identified modulo constants and represented by
zero-mean functions, it is natural to consider the space
$
C^\infty(\mathbb T^d)/\mathbb R\simeq C^\infty_\diamond(\mathbb T^d)
$.
Restriction identifies $H^{-1}_\diamond$ with
$(H^1_\diamond)^*$ (the dual of $H^1_\diamond$), with equivalent standard norms.

On $Q$, we use the parabolic H\"older convention
$C^{q,\beta}(Q):=C^{(q+\beta)/2,\,q+\beta}(Q)$,
with $C^\beta(Q):=C^{\beta/2,\,\beta}(Q)$.
The exponents on the right in the parabolic convention specify
temporal and spatial regularity, respectively; one time derivative
has weight two. For a finite-dimensional Euclidean space $E$,
these spaces are defined componentwise for $E$-valued functions.
All normed spaces carry their usual norms,
denoted by $\|\cdot\|_F$, with domains
omitted when clear.
We write $\mathcal L(E,F)$ for bounded linear operators, with
$\mathcal L(E):=\mathcal L(E,E)$, and use the operator norm. The notation
$A:E\xrightarrow{\ \sim\ }F$ means that
$A$ is a bounded linear isomorphism with
bounded inverse.
We use ${}^*$
for dual spaces or adjoints, as appropriate.
The symbol $D$ denotes differentiation in the indicated variables.
The notation $V\Subset W$ means that $\overline V$ is compact and
contained in $\operatorname{int}W$. Finally, whenever we write $\operatorname{int}U$ for a
control set $U\subseteq\mathbb E^d_+$,
the interior is understood with respect to
the full Euclidean space
$\mathbb E^d.$\\

We now introduce the semimartingale optimal transport problem
and its dual formulation, which provide the framework for
the mirror flow studied below.
Fix initial and terminal probability measures $\mu,\nu$ on
$\mathbb T^d$. 
Let $\Omega:=C([0,1];\mathbb T^d)$ and let
$X=(X_t)_{t\in[0,1]}$ be the canonical process.
For any probability law $P$ on $\Omega$, write
$P_t:=P\circ X_t^{-1}$.
Let $U\subseteq\mathbb E^d_+$ be a closed convex
set of admissible diffusion--drift pairs, with
nonempty interior. We denote by $\mathcal P(\mu)$
the set of probability laws $P$ on $\Omega$ such
that $P_0=\mu$ and $X$ admits the weak dynamics
\[
\mathrm dX_t
=
b_t^P\,\mathrm dt+\sqrt{2a_t^P}\,\mathrm dW_t,
\qquad
(a_t^P,b_t^P)\in U
\quad \mathrm dt\otimes P\text{-a.e.},
\]
where $a^P$ and $b^P$ are predictable processes
and $W$ is a standard $d$-dimensional Brownian
motion. The stochastic dynamics on the torus are
understood modulo $\mathbb Z^d$. The coefficient $a$ is half the instantaneous
covariance, and the generator for Markovian coefficients is
$b\cdot\nabla+a:\nabla^2$.
\noindent
For a Markovian running cost $L=L(t,x,a,b)$, define
\[
\mathcal C(P)
:=
\mathbb E^P\left[
\int_0^1 L(t,X_t,a_t^P,b_t^P)\,\mathrm dt
\right],
\qquad
V(\nu)
:=
\inf_{\substack{P\in\mathcal P(\mu)\\P_1=\nu}}
\mathcal C(P).
\]
For a terminal potential \(\psi\in C(\mathbb T^d)\), the penalized control value is
\begin{equation}
\label{eq:stoch_A}
\mathcal A(\psi)
:=
\inf_{P\in\mathcal P(\mu)}
\left\{
\mathcal C(P)+\int_{\mathbb T^d}\psi\,\mathrm dP_1
\right\}.
\end{equation}
Relaxing the terminal constraint with the multiplier $\psi$
gives the dual functional
\[
D_\nu(\psi)
:=
\inf_{P\in\mathcal P(\mu)}
\left\{
\mathcal C(P)+\int_{\mathbb T^d}\psi\,\mathrm d(P_1-\nu)
\right\}
=
\mathcal A(\psi)-\int_{\mathbb T^d}\psi\,\mathrm d\nu.
\]
The associated dual problem is then
\begin{equation}
\label{eq:dual}
\sup_{\psi\in C(\mathbb T^d)}D_\nu(\psi).
\end{equation}
\noindent
Define the Hamiltonian
\begin{equation}
\label{eq:H_control_representation}
H(t,x,M,p)
:=
\inf_{(a,b)\in U}
\left\{
a:M+b\cdot p+L(t,x,a,b)
\right\},
\qquad (M,p)\in\mathbb E^d,
\end{equation}
and denote its minimizer, when unique, by
$\mathfrak z=(\mathfrak a,\mathfrak b)$.
Under appropriate regularity assumptions, for a
sufficiently smooth terminal potential $\psi$, the value
function $u^\psi$ solves the Hamilton--Jacobi--Bellman equation
\begin{equation}
\label{eq:HJB_T_differentiability}
\partial_tu^\psi
+
H(t,x,\nabla^2u^\psi,\nabla u^\psi)
=0,
\qquad
u^\psi(1,\cdot)=\psi,
\end{equation}
and the stochastic control value satisfies
$
\mathcal A(\psi)
=
\int_{\mathbb T^d}u^\psi(0,x)\,\mu(\mathrm dx).
$
The optimal feedback
\[
(a^\psi,b^\psi)(t,x)
:=
\mathfrak z\bigl(
t,x,\nabla^2u^\psi(t,x),\nabla u^\psi(t,x)
\bigr)
\]
induces marginal densities satisfying the Fokker-Planck equation
\begin{equation}
\label{eq:FP_T_differentiability}
\partial_tm^\psi
=
-\partial_i(b_i^\psi m^\psi)
+\partial_{ij}(a_{ij}^\psi m^\psi),
\qquad
m^\psi(0,x)\,\mathrm dx=\mu(\mathrm dx).
\end{equation}
We write $T(\psi)$ for the optimal terminal law. In the PDE setting,
\[
T(\psi)(\mathrm dx)=m^\psi(1,x)\,\mathrm dx.
\]
When working with densities, we use the same symbol for a measure
and its density with respect to $\mathrm dx$. In particular, we
then write $T(\psi)=m^\psi(1,\cdot)$.
The terminal law $T(\psi)$ has unit mass, so its variations have zero mass.
Adding a constant to $\psi$ adds the same constant to
$u^\psi$, leaving the optimal feedback and its marginal
densities unchanged. Hence
$
T(\psi+c)=T(\psi),
c\in\mathbb R,
$
which justifies fixing the additive constant by working
with zero-mean terminal potentials.
Finally, under appropriate differentiability assumptions, the envelope identities are
\[
D_MH(t,x,M,p)=\mathfrak a(t,x,M,p),
\qquad
D_pH(t,x,M,p)=\mathfrak b(t,x,M,p).
\]

\medskip 

To compare terminal potentials and study the convergence of the numerical scheme proposed in this article for problem \eqref{eq:dual}, we need a measure of discrepancy adapted to the convex
structure of the problem. For this purpose, we introduce the
Bregman divergence.

\begin{definition}[Bregman divergence]
\label{def:bregman-divergence}
Let $\Phi$ be a convex, G\^ateaux differentiable functional
on a convex set $\mathcal D$ of zero-mean terminal potentials.
For $\psi,\bar\psi\in\mathcal D$, the Bregman divergence
generated by $\Phi$ is
\[
\calB_\Phi(\psi\mid\bar\psi)
:=
\Phi(\psi)-\Phi(\bar\psi)
-D\Phi(\bar\psi)[\psi-\bar\psi].
\]
\end{definition}

\noindent
The Bregman divergence generalizes the notion of squared distance
by measuring discrepancies through a convex functional $\Phi$. It measures the gap
between $\Phi(\psi)$ and its affine approximation at $\bar\psi$,
thereby adapting the measure of discrepancy to the curvature of
$\Phi$. It is nonnegative and, when $\Phi$ is strictly convex,
vanishes precisely when the two potentials coincide. However,
it need not be symmetric 
so it is generally not a metric. Its dependence on $\Phi$
allows the geometry of an optimization scheme to reflect the
structure of the underlying problem.

\section{Main Results}
\label{sec:main}
The goal of this paper is to introduce a mirror flow type scheme for solving \eqref{eq:dual} and to prove its convergence. 
We first establish the differentiability of the penalized value
$\mathcal A$ and use it to describe the Bregman geometry of the dual
problem~\eqref{eq:dual}. We then derive the associated mirror flow
and prove a contraction estimate for the exact discrete scheme.
Finally, we study the differentiability of the terminal map and
the invertibility of its derivative. These results allow us to
linearize the mirror step and bound the error of the open-loop
scheme uniformly in time.

\subsection{Construction of the mirror flow for semimartingale optimal transport}

Our goal is to construct a flow of terminal potentials whose
associated terminal marginals approach the target $\nu$ as fast as possible. To guide
this flow, we need to know how the dual objective changes when we
vary the potential. Since $D_\nu$ differs from $\mathcal A$ by an
affine term, the first step is to establish differentiability of
$\mathcal A$.

We begin by proving that the penalized problem attains its infimum.
This lets us compare its values at two potentials using their
respective optimal laws. The proofs use compactness and stability
arguments inspired by \cite{TaTo13}.

\begin{assumption}[Attainment and uniqueness conditions]
\label{ass:attainment_uniqueness}
Suppose the following hold.
\begin{enumerate}
\item[{\rm (i)}]
The running cost $L:[0,1]\times\mathbb T^d\times U\to[0,\infty)$
is continuous, convex in $(a,b)$, and uniformly continuous in $t$.
There exist $p>1$ and $C_0>0$ such that for $(t,x,a,b)\in[0,1]\times\mathbb T^d\times U$,
$
|a|^p+|b|^p
\leq C_0\bigl(1+L(t,x,a,b)\bigr)$.

\item[{\rm (ii)}]
The transport value $V$ is strictly convex on its effective domain.
\end{enumerate}
\end{assumption}

\begin{remark}[Finite-cost feasibility]
\label{rem:finite-cost-feasibility}
Assumption \ref{ass:attainment_uniqueness} implies the existence of an admissible law with finite running cost. Indeed, since $U$ has nonempty interior,
we may choose a constant pair
$(a_0,b_0)\in\operatorname{int}U$, with $a_0$
positive definite. The diffusion
$
dX_t=b_0\,dt+\sqrt{2a_0}\,dW_t$, 
$X_0\sim\mu,
$
on $\mathbb T^d$ defines an admissible law $P^0$.
By continuity of $L$ and compactness of
$[0,1]\times\mathbb T^d$,
$
\mathcal C(P^0)
\leq
\sup_{(t,x)\in[0,1]\times\mathbb T^d}
L(t,x,a_0,b_0)<\infty.
$
Thus the penalized problem always admits at least
one finite-cost competitor. 
\end{remark}

\begin{proposition}[Attainment of the penalized control problem]
\label{prop:attainment_A}
Under Assumption~\ref{ass:attainment_uniqueness}\textup{(i)},
for every $\psi\in C(\mathbb T^d)$, the value $\mathcal A(\psi)$
is finite and the infimum defining it is attained by some
$P^\psi\in\mathcal P(\mu)$.
\end{proposition}

\begin{proposition}[Uniqueness of the optimal terminal marginal]
\label{lem:terminal_marginal_uniqueness}
Let Assumption \ref{ass:attainment_uniqueness} hold. Then, for
every terminal potential $\psi$, all minimizers of $\mathcal A(\psi)$ have the
same terminal marginal. We denote this marginal by $T(\psi)$.
\end{proposition}

\begin{proposition}[Differentiability of the penalized control value]
\label{prop:A_diff}
Under Assumption~\ref{ass:attainment_uniqueness},
the functional
$
\mathcal A:C(\mathbb T^d)\longrightarrow\mathbb R
$
is G\^ateaux differentiable. For every
$\psi,h\in C(\mathbb T^d)$,
\[
D\mathcal A(\psi)[h]
=
\langle T(\psi),h\rangle,
\]
where $T(\psi)$ is the unique optimal terminal
marginal given by Proposition \ref{lem:terminal_marginal_uniqueness}.
\end{proposition}

\begin{remark}[Stochastic and PDE representations of the terminal map]
Under Assumptions~\ref{ass:attainment_uniqueness}
and~\ref{ass:T_diff}, for every
$\psi\in C^{m+2,\alpha}_\diamond(\mathbb T^d)$,
the HJB--Fokker--Planck system
\eqref{eq:HJB_T_differentiability}--%
\eqref{eq:FP_T_differentiability}
admits the classical solution described in
Lemma~\ref{lem:state_system_terminal_map}.
The HJB verification argument shows that the
diffusion with feedback coefficients
$(a^\psi,b^\psi)$ minimizes
\eqref{eq:stoch_A}.
By Proposition~\ref{lem:terminal_marginal_uniqueness},
its terminal law is therefore $T(\psi)$.
Consequently,
\[
T(\psi)(\mathrm dx)=m^\psi(1,x)\,\mathrm dx
\]
with $m^{\psi}$ given by \eqref{eq:FP_T_differentiability}.
In this regular setting, we also use
$T(\psi)$ to denote its density
$m^\psi(1,\cdot)$.
\end{remark}

\medskip

\noindent
Proposition \ref{prop:A_diff} directly yields
$
D(D_\nu)(\psi)[h]
=
\pair{T(\psi)-\nu}{h}.
$
The first variation thus identifies the terminal mismatch
$T(\psi)-\nu$. Since $\mathcal A$ is the pointwise infimum of affine
functionals of $\psi$, it is concave.
To measure changes in the potential, we apply Definition~\ref{def:bregman-divergence} with $\Phi:=-\mathcal A$ and obtain
\[
\calB_\Phi(\psi\mid\bar\psi)
=
-\mathcal A(\psi)+\mathcal A(\bar\psi)
+\pair{T(\bar\psi)}{\psi-\bar\psi}.
\]
Since $-D_\nu$ differs from $\Phi$ only by an affine term, it
generates the same divergence. The resulting geometry is therefore
independent of the target $\nu$.
In the differentiability regime studied further in this article, the operator
$-DT_\psi$ represents the Bregman Hessian:
$
D^2\Phi(\psi)[h,k]=\pair{(-DT_\psi) h}{k}.
$
Consequently, as $h\to0$ in $C^{m+2,\alpha}$,
\[
\calB_\Phi(\psi+h\mid\psi)
=
\frac12\pair{(-DT_\psi) h}{h}
+
o\!\left(\|h\|_{C^{m+2,\alpha}}^2\right).
\]
The negative Hessian of $\mathcal A$ thus determines the local
geometry of the potential space. The size assigned to a perturbation
$h$ at each step of the scheme depends on the curvature of the dual objective in that direction
at $\psi$, rather than on a fixed Euclidean norm.\\

Let's now introduce the continuous-time mirror flow and study the dual gap along it. Where $(-DT_\psi)$ is invertible on the zero-mean spaces, we define
the continuous mirror flow by
\[
\frac{\mathrm d}{\mathrm ds}D\Phi(\psi_s)
=
D(D_\nu)(\psi_s).
\]
Let $p_s:=T(\psi_s)$ denote the terminal marginal.
Since $D\Phi(\psi_s)=-p_s$ and
$D(D_\nu)(\psi_s)=p_s-\nu$, the flow gives
$
\frac{\mathrm dp_s}{\mathrm ds}=\nu-p_s.
$
Thus the terminal marginal converges exponentially fast toward the target $\nu$
\[
p_s=\nu+e^{-s}(p_0-\nu),
\qquad p_0:=T(\psi_0).
\]
\noindent
Suppose that a potential $\psi^\star$ satisfies $T(\psi^\star)=\nu$.
Then $\psi^\star$ maximizes the concave objective $D_\nu$, and
the Bregman divergence with reference potential $\psi^\star$
equals the dual gap:
$
\calB_\Phi(\psi\mid\psi^\star)
=
D_\nu(\psi^\star)-D_\nu(\psi)
$. This identity makes the Bregman divergence a natural choice of
distance for our algorithm: moving closer to $\psi^\star$ in this
geometry is exactly equivalent to reducing the dual gap. We can
therefore use this divergence to track progress toward a potential
that produces the target marginal $\nu$.
Finally, differentiating this identity along the flow gives
\[
\frac{\dd}{\dd s}\calB_\Phi(\psi_s\mid\psi^\star)
=
-\pair{p_s-\nu}{(-DT_{\psi_s})^{-1}(p_s-\nu)}
\leq 0,
\]
where the inequality follows from the nonnegativity of the
Bregman Hessian. The dual gap therefore serves as a Lyapunov
functional: it is nonincreasing along the flow and vanishes
whenever the terminal marginal matches $\nu$.

\subsection{Discretized Mirror Flow}
\noindent
We now discretize the continuous mirror flow to obtain an iterative
scheme for the terminal potential. For a step size $\tau>0$, we
consider a Bregman minimizing movement: each update balances
improvement of the dual objective with a penalty on the change in
potential, measured by the Bregman divergence. We work with smooth terminal potentials normalized to have zero mean.
Given $\psi_0\in\calX$, we denote the exact discrete iterates by
$\bar\psi_k$, with $\bar\psi_0:=\psi_0$, and reserve $\psi_k$ for the
linearized schemes introduced below. For each $k\geq0$, we set
\begin{equation}
\label{eq:bregman_mm}
\bar\psi_{k+1}
\in
\argmin_{\psi\in\calX}
\left\{
-D_\nu(\psi)
+\frac{1}{\tau}\calB_\Phi(\psi\mid\bar\psi_k)
\right\}.
\end{equation}
This scheme corresponds to an implicit time discretization of the
mirror flow in the coordinates $D\Phi(\psi)$.
Using
$
D\Phi(\psi)=-T(\psi)
$
and
$
D(D_\nu)(\psi)=T(\psi)-\nu,
$
the first-order optimality condition yields
$
T(\bar\psi_{k+1})-T(\bar\psi_k)
=
\tau\bigl(\nu-T(\bar\psi_{k+1})\bigr),
$
and therefore
\[
T(\bar\psi_{k+1})
=
\frac{1}{1+\tau}T(\bar\psi_k)
+
\frac{\tau}{1+\tau}\nu.
\]

\noindent
Motivated by this recurrence, we define the prescribed
terminal-marginal path independently of the existence of the exact
mirror iterates by
\begin{equation}
\label{eq:explicit_flow}
\bar p_k
:=
\nu+(1+\tau)^{-k}(\bar p_0-\nu),
\qquad
\bar p_0:=T(\psi_0).
\end{equation}
Whenever the exact mirror iterates $\bar\psi_k$ exist, their
terminal marginals satisfy $T(\bar\psi_k)=\bar p_k$.

The prescribed path is determined by the initial terminal marginal
$\bar p_0$, the target marginal $\nu$, and the step size $\tau$.
It therefore defines a prescribed
open-loop path in terminal coordinates, along which the terminal
mismatch decreases geometrically with contraction factor
$(1+\tau)^{-1}$. We now quantify progress
toward an optimal potential in the Bregman geometry and prove the same contraction for the reverse Bregman divergence
$\calB_\Phi(\psi^\star\mid\bar\psi_k)$.

\begin{theorem}[Discrete reverse-Bregman contraction]
\label{thm:Bregman_contraction}
Suppose that Assumption~\ref{ass:attainment_uniqueness}
holds. Fix $\tau>0$ and an initial potential
$\bar\psi_0\in\calX$. Assume that, for every $k\geq0$, the minimization
problem \eqref{eq:bregman_mm} admits a minimizer
$\bar\psi_{k+1}\in\calX$, and that there exists
$\psi^\star\in\calX$ satisfying
$T(\psi^\star)=\nu$. Then, for every $k\geq0$,
\[
\calB_\Phi(\psi^\star\mid\bar\psi_k)
\leq
(1+\tau)^{-k}
\calB_\Phi(\psi^\star\mid\bar\psi_0).
\]
\end{theorem}

\noindent
We also want to know how close the current dual value is to
the maximum. We establish that the dual
objective gap converges geometrically to zero and give an
explicit bound on the remaining optimization error after
$k$ iterations.

\begin{corollary}[Convergence of the objective value]
\label{cor:value_convergence}
Under the assumptions of the previous theorem, for every $k\geq1$,
\begin{equation}
\label{eq:value_convergence}
D_\nu(\psi^\star)-D_\nu(\bar\psi_k)
\leq
\frac{1}{\tau}(1+\tau)^{-(k-1)}
\calB_\Phi(\psi^\star\mid\bar\psi_0).
\end{equation}
\end{corollary}

\begin{remark}[Convergence of $\Phi$]
Since
$
\Phi(\bar\psi_k)-\Phi(\psi^\star)
=
D_\nu(\psi^\star)-D_\nu(\bar\psi_k)
-
\pair{\nu}{\bar\psi_k-\psi^\star},
$
convergence of $D_\nu(\bar\psi_k)$ implies convergence of
$\Phi(\bar\psi_k)$ provided
$
\pair{\nu}{\bar\psi_k-\psi^\star}\to0.
$
For example, this follows if there exist a norm $\|\cdot\|_Z$
and a constant $c>0$ such that
$\calB_\Phi(\psi^\star\mid\bar\psi_k)
\geq c\|\bar\psi_k-\psi^\star\|_Z^2$ for every $k$,
and $\psi\mapsto\langle\nu,\psi\rangle$ is continuous in this norm.
\end{remark}

\subsection{Linearized mirror step}
To implement the exact scheme \eqref{eq:explicit_flow}, it remains to recover the potential
associated with each prescribed terminal marginal.
If the nonlinear inverse of $T$ could be computed easily, the prescribed path
could be lifted exactly to potential space by setting
$
\bar\psi_k:=T^{-1}(\bar p_k).
$
Equivalently, one would solve the nonlinear equation
$
T(\bar\psi_{k+1})=\bar p_{k+1}
$
at every step. We instead consider two linearized approximations of
this exact lift. We first describe the two updates formally, assuming that $T$ is
differentiable and that $DT_\psi$ is invertible on the relevant
zero-mean spaces. These properties will be established further.

\paragraph{Prescribed open-loop linearization.}
We seek to reproduce, to first order, the successive changes
along the prescribed reference path $(\bar p_k)$. Define
$
\rho_0:=\bar p_0-\nu$ and 
$d_k
:=
\bar p_k-\bar p_{k+1}
=
\tau(1+\tau)^{-k-1}\rho_0$.
Writing $\psi_{k+1}=\psi_k+h_k$, we choose the potential increment
so that the linearized change in the actual terminal density
matches the prescribed change:
$
DT_{\psi_k}[h_k]
=
\bar p_{k+1}-\bar p_k
=
-d_k.
$
This gives the open-loop recursion
\begin{equation}
\label{eq:recursion}
h_k^{\mathrm{ol}}
:=
-(DT_{\psi_k})^{-1}d_k,
\qquad
\psi_{k+1}
:=
\psi_k+h_k^{\mathrm{ol}}.
\end{equation}
The terminal-coordinate increments $d_k$ are fixed in advance,
while the inverse $DT_{\psi_k}^{-1}$ is evaluated at the current
potential. In particular, the prescribed increments are not
adjusted according to the difference between $T(\psi_k)$ and
$\bar p_k$. Tracking errors accumulated during earlier steps
are therefore not corrected by this update.

\paragraph{Adaptive linearization.}
Alternatively, we may choose the update using the terminal
density actually reached at step $k$. We ask the linearized
terminal density to correct the fraction $\tau/(1+\tau)$ of
the current mismatch with the target:
$
T(\psi_k)+DT_{\psi_k}[h_k]
=
T(\psi_k)
+
\frac{\tau}{1+\tau}\bigl(\nu-T(\psi_k)\bigr).
$
This yields the adaptive recursion
\begin{equation}
\label{eq:adaptive-linearized-recursion}
h_k^{\mathrm{ad}}
:=
\frac{\tau}{1+\tau}
(DT_{\psi_k})^{-1}
\bigl(\nu-T(\psi_k)\bigr),
\qquad
\psi_{k+1}
:=
\psi_k+h_k^{\mathrm{ad}}.
\end{equation}
Here the direction in terminal coordinates is recomputed from
the current mismatch, allowing the update to respond to the
actual evolution of the terminal density.\\

To compare the two update rules at a common current potential
$\psi_k$, set
$
\theta_k:=T(\psi_k)-\bar p_k$. The adaptive direction in terminal coordinates is
$
d_k^{\mathrm{ad}}
:=
\frac{\tau}{1+\tau}\bigl(T(\psi_k)-\nu\bigr).
$
Since
$
T(\psi_k)-\nu
=
\bar p_k-\nu+\theta_k,
$
we obtain
$
d_k^{\mathrm{ad}}
=
d_k+\frac{\tau}{1+\tau}\theta_k.
$
Consequently,
\[
h_k^{\mathrm{ol}}=(-DT_{\psi_k})^{-1}d_k,
\qquad
h_k^{\mathrm{ad}}
=
(-DT_{\psi_k})^{-1}
\left(
d_k+\frac{\tau}{1+\tau}\theta_k
\right).
\]
Thus the adaptive rule adds a feedback correction proportional
to the deviation from the prescribed reference path. Starting
both schemes from the same potential $\psi_0$, their first
updates coincide because $T(\psi_0)=\bar p_0$, hence $\theta_0=0$.
Their subsequent trajectories may differ.

\vspace{\baselineskip}

We now justify the first-order approximation used in the
linearized schemes. We first establish differentiability of
the terminal map $T$ and estimate the error of its linear
approximation. Under
additional smoothness assumptions, we then prove that
$DT_\psi$ is an isomorphism at smooth terminal potentials.
Finally, we show that smoothness propagates along the
linearized recursion, ensuring that every iteration is
well defined.

For this analysis, we restrict to an autonomous,
state-independent running cost $L=L(a,b)$. The Hamiltonian
introduced in \eqref{eq:H_control_representation} becomes
\begin{equation}
\label{eq:autonomous_H_control_representation}
H(M,p)
:=
\inf_{(a,b)\in U}
\left\{
a:M+b\cdot p+L(a,b)
\right\},
\qquad (M,p)\in\mathbb E^d.
\end{equation}
We denote its optimizer by
$\mathfrak z(M,p)
=(\mathfrak a(M,p),\mathfrak b(M,p))$. Accordingly, throughout this autonomous subsection,
\[
\partial_tu^\psi+H(\nabla^2u^\psi,\nabla u^\psi)=0,
\qquad
(a^\psi,b^\psi)=\mathfrak z(\nabla^2u^\psi,\nabla u^\psi).
\]
This restriction isolates the main arguments from the additional
technicalities caused by explicit dependence on \(t\) and \(x\).  We expect the differentiability and invertibility theory developed below to extend to the general case \(L=L(t,x,a,b)\) under suitable regularity assumptions; the extension is mainly technical and is left for future work.

\begin{assumption}[Semimartingale regularity and ellipticity]
\label{ass:T_diff}
Assume the following.
\begin{enumerate}
\item[\rm(i)]
$L\in C^{m+3,\alpha}_{\mathrm{loc}}(\operatorname{int}U)$ as a function of
the full joint control $(a,b)$.

\item[\rm(ii)]
For every open set $\mathcal P\Subset\mathbb E^d$, there exists a
compact set $\mathcal K_{\mathcal P}\Subset U$ such that, for
every $(M,p)\in\mathcal P$, the infimum in
\eqref{eq:autonomous_H_control_representation} has a unique minimizer
$\mathfrak z(M,p)\in\mathcal K_{\mathcal P}$.

\item[\rm(iii)]
For every $\mathcal P$ as in \textup{(ii)}, there is a constant
$0<c_{L,\mathcal P}$ such that
\[
c_{L,\mathcal P}|(A,\beta)|^2
\leq D^2L(a,b)[(A,\beta),(A,\beta)]
\]
for all $(a,b)\in\mathcal K_{\mathcal P}$ and
$(A,\beta)\in\mathbb E^d$.

\item[\rm(iv)]
There are constants $0<\lambda\leq\Lambda<\infty$ such that
\[
\lambda I_d\leq D_MH(M,p)\leq\Lambda I_d,
\qquad (M,p)\in\mathbb E^d.
\]

\item[\rm(v)]
The initial probability measure $\mu$ admits a density
in $C^{m,\alpha}(\mathbb T^d)$.
\end{enumerate}
\end{assumption}

\begin{remark}[Role of the assumptions]
\begin{enumerate}
\item We require the optimizer to lie in $\operatorname{int}U$
so that its optimality condition can be written as
$
DL\bigl(\mathfrak z(M,p)\bigr)+(M,p)=0.
$
Differentiating this equation then allows us to determine
how the optimal controls vary with the terminal potential.
When the gradient and Hessian of the value function remain
uniformly bounded, the associated optimal controls lie in
a common compact subset of $\operatorname{int}U$. This compact localization provides uniform bounds on the derivatives
of $L$ and, together with positive definiteness, on the inverse
Hessian. It also keeps the controls a positive distance from
the boundary. These bounds are used in the implicit-function
and Taylor arguments and in obtaining uniform local estimates
for the associated parabolic equations.
Assumption~\ref{ass:T_diff} collects sufficient conditions
for this construction.

\item Under assumptions (i)--(iii), the interior
optimality condition and the implicit function
theorem imply
$
\mathfrak z\in C^{m+2,\alpha}_{\mathrm{loc}}
(\mathbb E^d).
$
Since $DH=\mathfrak z$, it follows that
$
H\in C^{m+3,\alpha}_{\mathrm{loc}}
(\mathbb E^d).
$
Moreover, the upper bounds on $D^2L$ over
compact optimizer sets follow from the
local regularity of $L$ and compactness. These regularity properties are established
explicitly at the beginning of the proof of
Lemma~\ref{lem:state_system_terminal_map}.
\end{enumerate}
\end{remark}

\begin{theorem}[Strong differentiability of the terminal map]
\label{prop:frechet_terminal_map}
Under Assumption~\ref{ass:T_diff}, the terminal map
$
T:C^{m+2,\alpha}_\diamond(\mathbb T^d)
\longrightarrow C^{m,\alpha}(\mathbb T^d)
$
is locally of class $C^{1,1}$. 
For every
$\psi,h\in C^{m+2,\alpha}_\diamond(\mathbb T^d)$,
the linearized system
\eqref{eq:linearized_system} admits a unique
solution
\[
(w^h,\dot z^h,q^h)\in
C^{m+2,\alpha}(Q)
\times C^{m,\alpha}(Q;\mathbb E^d)
\times C^{m,\alpha}(Q),
\qquad
\dot z^h=(\dot a^h,\dot b^h).
\]
Moreover, the following representation holds:
\[
DT_\psi[h]=q^h(1,\cdot)
\in C^{m,\alpha}_\diamond(\mathbb T^d).
\]
\end{theorem}

\begin{theorem}[Hölder invertibility of the linearized terminal map at smooth potentials]
\label{prop:strong-inverse-smooth-terminal}
Under Assumption~\ref{ass:T_diff}, suppose that
$L\in C^\infty_{\mathrm{loc}}(\operatorname{int}U)$.
Then, for every
$\psi\in C^\infty_\diamond(\mathbb T^d)$,
$
DT_\psi:
C^{m+2,\alpha}_\diamond(\mathbb T^d)
\xrightarrow{\ \sim\ }
C^{m,\alpha}_\diamond(\mathbb T^d)
$
is an isomorphism. Moreover, its inverse maps
$C^\infty_\diamond(\mathbb T^d)$
into $C^\infty_\diamond(\mathbb T^d)$.
\end{theorem}

\noindent
For the open-loop recursion \eqref{eq:recursion}, we write
$h_k:=h_k^{\mathrm{ol}}$ in the results and proofs below. We first show
that each iteration is well defined.

\begin{corollary}
[Well-definedness and smoothness of the linearized recursion]
\label{prop:smooth-induction-geometric-directions}
Under Assumption~\ref{ass:T_diff}, suppose that
$L\in C^\infty_{\mathrm{loc}}(\operatorname{int}U)$ and
$\psi_0,\rho_0\in C^\infty_\diamond$.
Then the recursion \eqref{eq:recursion} is well defined for
every $k\geq0$, with
$
\psi_k,h_k\in C^\infty_\diamond.
$
\end{corollary}

\subsection{Uniform-in-time error of the open-loop linearization}

We now quantify the error produced by replacing the nonlinear inverse of the terminal map by the linearized
recursion. We first bound the error introduced by a single
update, then control its accumulation over successive iterations.
Under local regularity and smallness assumptions around $\psi_0$,
we show that the tracking error remains bounded in $C^{m,\alpha}$
uniformly over all iterations and is of order $\tau$ as
$\tau\downarrow0$.

\begin{assumption}[Additional conditions for the error analysis]
\label{ass:open-loop-error}
Let Assumption~\ref{ass:T_diff} hold.
In addition,
assume that
\[
L\in C^\infty_{\mathrm{loc}}(\operatorname{int}U),
\qquad
\psi_0\in C^\infty_\diamond(\mathbb T^d),
\]
and that the target probability measure $\nu$ admits
a density in $C^{m,\alpha}(\mathbb T^d)$.
\end{assumption}
\begin{remark}
Under Assumption \ref{ass:open-loop-error} setting
$
\rho_0:=T(\psi_0)-\nu,
$
we have
$\rho_0\in C^{m,\alpha}_\diamond(\mathbb T^d)$.
\end{remark}







\noindent
We work in a neighborhood of the reference potential $\psi_0$,
written
\[
\calK_R
:=
\left\{
\psi\in C^{m+2,\alpha}_\diamond(\mathbb T^d):
\|\psi-\psi_0\|_{C^{m+2,\alpha}}\leq R
\right\}.
\]
For sufficiently small $R$, the assumptions yield uniform
local estimates on this ball.


\begin{theorem}[Uniform-in-time open-loop tracking error]
\label{thm:open_loop_error}
Let Assumption~\ref{ass:open-loop-error} hold. 
There exist $R_0,M, C_{\mathrm{err}}>0$ with the following property.
For every
$R\in(0,R_0]$, if $\tau>0$ and
$
M\|\rho_0\|_{C^{m,\alpha}}\leq R,
$
then the recursion \eqref{eq:recursion} is well defined
for every $k\geq0$, its iterates remain in $\calK_R$, and
\[
\sup_{k\geq0}
\|(-DT_{\psi_k})^{-1}\|_{\mathcal L(
C^{m,\alpha}_\diamond,
C^{m+2,\alpha}_\diamond)}\leq M.
\]
Moreover,
\begin{equation}
\label{eq:holder-uniform-tracking}
\sup_{N\geq0}
\|T(\psi_N)-\bar p_N\|_{C^{m,\alpha}}
\leq C_{\mathrm{err}}
\frac{\tau}{2+\tau}
\|\rho_0\|_{C^{m,\alpha}}^2.
\end{equation}
Consequently, for every $N\geq0$,
\begin{equation}
\label{eq:holder-target-error}
\|T(\psi_N)-\nu\|_{C^{m,\alpha}}
\leq
(1+\tau)^{-N}\|\rho_0\|_{C^{m,\alpha}}
+C_{\mathrm{err}}\frac{\tau}{2+\tau}
\|\rho_0\|_{C^{m,\alpha}}^2.
\end{equation}
The constants $R_0$, $M$, and $C_{\mathrm{err}}$
are independent of $\tau$, $\rho_0$, and
the iteration index $N$.
\end{theorem}

\section{Numerical experiments}
\label{sec:numerical-experiments}

We consider the one-dimensional periodic SEMOT model of
\cite[Section~4.2]{BuetGolfouseEtAl2026SEMOT}, with zero drift and a
Brownian reference of variance one. Our goal is to use the adaptive mirror scheme \eqref{eq:adaptive-linearized-recursion} to solve the
SEMOT dual problem and compare its convergence with ordinary gradient
flow. Both methods start from the same terminal potential and aim to
recover the target density $\nu$. By tracking the terminal mismatch
and the dual objective over the iterations, we illustrate how
adapting the updates to the local geometry accelerates convergence.
\noindent
We work in the Markovian setting on
$[0,t_{\mathrm f}]\times\mathbb T$, where $\mathbb T=\mathbb R/\mathbb Z$
and $t_{\mathrm f}=0.1$, and write the problem as
\[
\inf_{\mathbb P\in\mathcal P(\mu,\nu)}
\mathbb E^{\mathbb P}\!\left[
\int_0^{t_{\mathrm f}}\ell(\sigma_t^{\mathbb P})\,\mathrm dt
\right],
\qquad
\ell(\sigma):=\sigma\log\sigma-\sigma+1,\quad \ell(0)=1.
\]
Here $\mathcal P(\mu,\nu)$ denotes the admissible laws of driftless
Markovian diffusions on $\Omega_{\mathrm f}:=C([0,t_{\mathrm f}];\mathbb T)$ with initial marginal $\mu$ and
terminal marginal $\nu$ and
$\sigma$ is the instantaneous variance, so $a=\sigma/2$ in the
notation of the preceding sections. The Hamiltonian is
$H(M)=1-e^{-M/2}$. For a terminal potential $\psi$, the optimal
variance and marginal density satisfy
\begin{equation}
\begin{cases}
\partial_tu^\psi+1-e^{-\partial_{xx}u^\psi/2}=0,
&u^\psi(t_{\mathrm f})=\psi,\\
\sigma^\psi=e^{-\partial_{xx}u^\psi/2},&\\
\partial_tm^\psi=\tfrac12\partial_{xx}(\sigma^\psi m^\psi),
&m^\psi(0)=\mu.
\end{cases}
\label{eq:num-state-system}
\end{equation}
Writing $T(\psi)=m^\psi(t_{\mathrm f})$, the dual objective is
\[
D_{\nu}(\psi)
=\int_{\mathbb T}u^\psi(0,x)\mu(x)\,\mathrm dx
-\int_{\mathbb T}\psi(x)\nu(x)\,\mathrm dx,
\qquad
D (D_{\nu})(\psi)=T(\psi)-\nu.
\]
\noindent
The analytical assumptions of the preceding
sections may be viewed as describing a
uniformly elliptic, smooth-control
regularization of SEMOT. Here, however,
we retain the original entropy cost,
without imposing such a regularization.
For the
model considered in this section, the corresponding local regularity,
invertibility, and open-loop error estimates
can instead be established near the reference
solution using a model-specific curvature
maximum principle. These additional
arguments are not included in the general
analysis presented here.
\noindent
At a current potential $\psi_k$, set $r_k=T(\psi_k)-\nu$. Ordinary gradient ascent uses
\[
\psi_{k+1}=\psi_k+\eta_{\mathrm{GF}}r_k,
\]
whereas the adaptive mirror scheme uses
\[
\psi_{k+1}=\psi_k+h_k,\qquad
(-DT_{\psi_k})h_k=\eta_{\mathrm{MF}}r_k.
\]
Here $\eta_{\mathrm{GF}}$ and $\eta_{\mathrm{MF}}$ are the respective
outer step sizes. Setting $\eta_{\mathrm{MF}}=\tau/(1+\tau)$ recovers
\eqref{eq:adaptive-linearized-recursion}.
Both methods use the terminal mismatch $r_k$, which is the gradient
of the dual objective. Gradient flow takes a step in this direction,
while the mirror scheme applies $(-DT_{\psi_k})^{-1}$ to the gradient,
adjusting the update to the local curvature of the objective.
Writing $F=-D_{\nu}$, we have $\nabla F(\psi_k)=-r_k$ and
$\nabla^2F(\psi_k)=(-DT_{\psi_k})$, so the mirror increment becomes
\[
h_k=-\eta_{\mathrm{MF}}
\bigl(\nabla^2F(\psi_k)\bigr)^{-1}\nabla F(\psi_k).
\]
Thus the adaptive mirror scheme is a \emph{preconditioned gradient
descent} for $F$, with the inverse Hessian as preconditioner. For the numerical discretization, let $h=(\Delta x,\Delta t)$ denote
the mesh sizes. The spatial grid is $x_i=i\Delta x$,
$i=0,\ldots,N_x-1$, with $\Delta x=1/N_x$, and
$\Delta t=t_{\mathrm f}/N_t$. For a grid potential $\psi$,
$T_h(\psi)_i$ is the terminal density at $x_i$ obtained by solving
the discretized HJB--Fokker--Planck system. Let $\mu_h$ and
$\nu_h=(\nu_i)_{i=0}^{N_x-1}$ denote the initial and target densities
on the grid. Let $r_k:=T_h(\psi_k)-\nu_h$ be the discrete terminal mismatch,
with components $r_{k,i}$.

\begin{algorithm}[htbp]
\caption{Adaptive mirror scheme for SEMOT}
\label{alg:num-adaptive-mirror}
\begin{algorithmic}[1]
\Require $\mu_h,\nu_h$, initial potential $\psi_0=0$,
step size $\eta_{\mathrm{MF}}\in(0,1)$, tolerance $\varepsilon$
\For{$k=0,1,\ldots$}
    \State Solve the discretized system~\eqref{eq:num-state-system}
    and set $r_k=T_h(\psi_k)-\nu_h$
    \If{$\Delta x\sum_{i=0}^{N_x-1}|r_{k,i}|\leq\varepsilon$}
        \State \Return $\psi_k$
    \EndIf
    \State Assemble $DT_h(\psi_k)$ from the linearized discrete system
    \State Compute $h_k$ as a minimizer of
    $\|DT_h(\psi_k)[\delta\psi]+\eta_{\mathrm{MF}}r_k\|_{\ell^2}^2$
    over $\delta\psi\in\mathbb R^{N_x}$ with
    $\sum_{i=0}^{N_x-1}\delta\psi_i=0$
    \State Set $\psi_{k+1}=\psi_k+h_k$ and subtract its grid mean
\EndFor
\end{algorithmic}
\end{algorithm}

\begin{figure}[H]
\centering
\includegraphics[width=\textwidth]{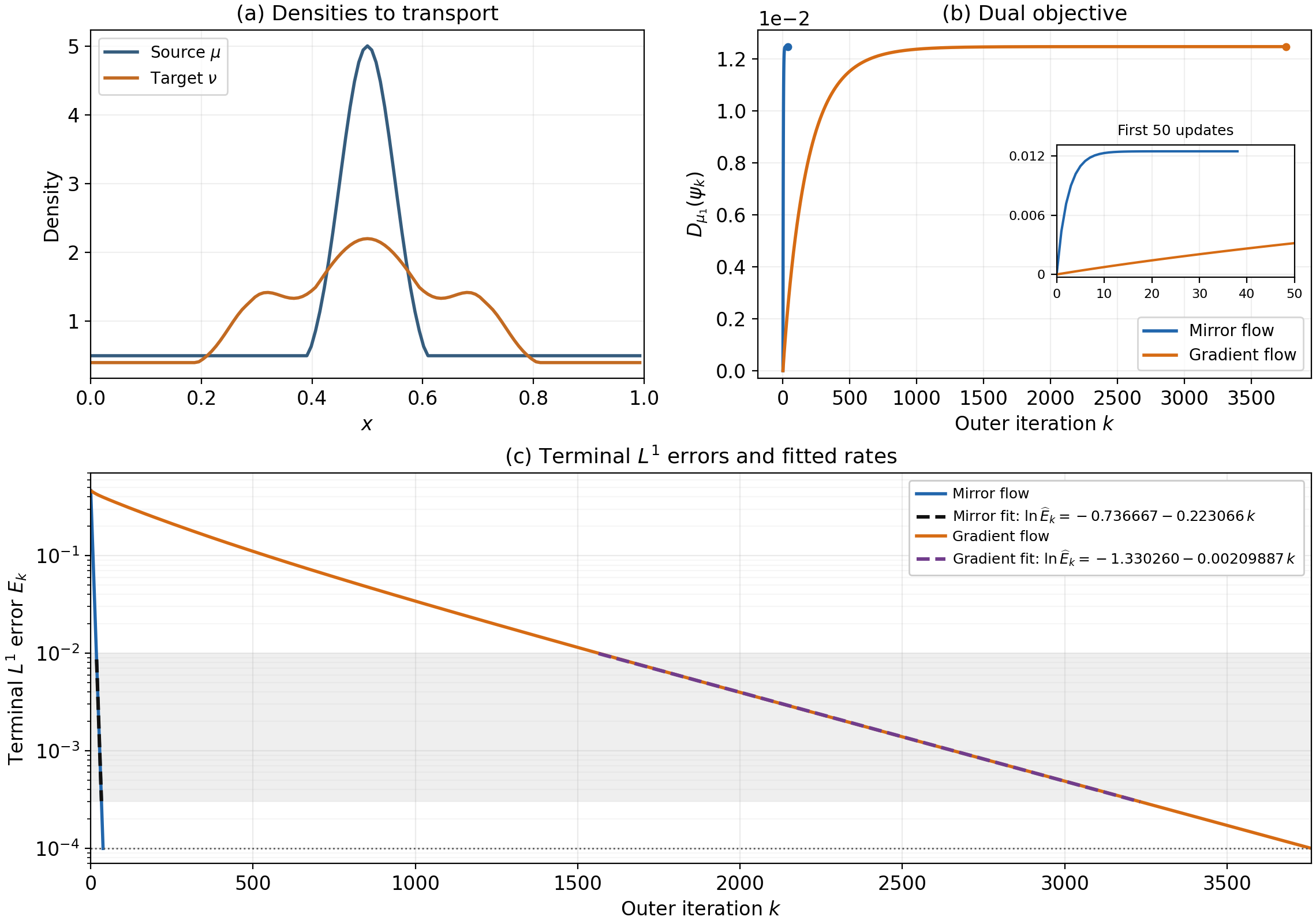}
\caption{Source and target densities (top left), dual objective
against iterations (top right, with the first 50 updates in the
inset), and both terminal $L^1$ error curves on common iteration
axes (bottom). Dashed lines show the exponentials
$E_k\simeq e^{b-\widehat r k}$ fitted over the grey shaded error range
$3\times10^{-4}\leq E_k\leq10^{-2}$.}
\label{fig:num-convergence}
\end{figure}

\medskip
\noindent
Algorithm~\ref{alg:num-adaptive-mirror} summarizes the mirror iteration.
At each iteration, both methods recompute the HJB and Fokker--Planck
solutions at the current potential and form the terminal mismatch
$r_k$. Gradient flow updates directly along this mismatch. The adaptive mirror scheme assembles $DT_h(\psi_k)$ from the
linearized discrete system and chooses a zero-mean increment $h_k$
that best fits $DT_h(\psi_k)[h_k]=-\eta_{\mathrm{MF}}r_k$
in the grid $\ell^2$ norm.
Both methods use the same initial potential, discretization, and
stopping criterion.
We use centered periodic finite differences and implicit time
stepping, with the forward scheme chosen as the discrete adjoint
of the HJB linearization. The mirror method uses least squares
on the zero-mean space. We use $N_x=128$ spatial nodes and $N_t=80$ time steps, so
$\Delta x=1/128$ and $\Delta t=0.1/80$.
The marginal densities are
\[
\mu=G_{0.5,0.05},
\qquad
\nu
=0.6G_{0.5,0.1}
+0.2G_{0.3,0.05}
+0.2G_{0.7,0.05},
\]
where $G_{m,s}$ denotes the sampled periodized Gaussian profile
centered at $m$ with width $s$, floored at $0.1$ before normalization.
Both methods start from $\psi_0=0$ and stop when
\[
E_k
:=\Delta x\sum_i
\left|T_h(\psi_k)_i-\nu_{i}\right|
\leq10^{-4}.
\]
We take $\eta_{\mathrm{MF}}=0.2$ for mirror flow and
$\eta_{\mathrm{GF}}=2.96868\times10^{-4}$ for gradient ascent.
\noindent
We estimate the exponential convergence rate $\widehat r$ by
ordinary least-squares fitting of
\[
\log E_k\simeq b-\widehat r\,k,
\qquad
E_k\simeq e^b\widehat q^{\,k},
\qquad
\widehat q=e^{-\widehat r},
\]
over the common error range
$3\times10^{-4}\leq E_k\leq10^{-2}$.
Here $\log$ denotes the natural logarithm, so the fitted slope is
$-\widehat r$. The fitted lines are
\begin{align*}
\log\widehat E_k^{\mathrm{MF}}
&=-0.736667-0.223066\,k,
& R^2&=0.99999998,\\
\log\widehat E_k^{\mathrm{GF}}
&=-1.330260-0.00209887\,k,
& R^2&=0.99999426.
\end{align*}
Both fits are excellent within this error range.
A larger $\widehat r$, or equivalently a smaller contraction factor
$\widehat q$, indicates faster convergence. Table~\ref{tab:num-comparison} reports the update counts, runtimes,
fitted rates, and final errors. Figure~\ref{fig:num-convergence} shows
the corresponding error trajectories and fitted exponential regimes.

\begin{table}[htbp]
\centering
\resizebox{\textwidth}{!}{%
\begin{tabular}{lrrrrrr}
\hline
Method & Updates & Time (s) & Rate $\widehat r$
& Factor $\widehat q$ & Terminal error & Dual value\\
\hline
Mirror flow & 38 & 3.26 & 0.22307 & 0.80006
& $9.9662\times10^{-5}$ & 0.0124614549\\
Gradient flow & 3,759 & 115.31 & 0.002099 & 0.99790
& $9.9998\times10^{-5}$ & 0.0124614547\\
\hline
\end{tabular}%
}
\caption{Convergence to $E_k\leq10^{-4}$.
Rates and contraction factors are estimated per outer update.
Terminal errors $E_N$ and dual values $D_{\nu}(\psi_N)$
are evaluated at the final iterate.
Times are from one local run and include all state solves and updates.}
\label{tab:num-comparison}
\end{table}

Mirror flow reaches the tolerance with about 99 times fewer updates
and a measured 35-fold speedup in this run.
Its contraction factor $0.80006$ agrees closely with the predicted
$1-\eta_{\mathrm{MF}}=0.8$, whereas the gradient factor $0.99790$
is close to one. These results support the local spectral explanation
of the benefit of preconditioning.
Indeed, in this case the numerical difficulty is the poor conditioning of the discrete
$L^2$ Hessian which has a large condition number. On the zero-mean grid space, the discrete Hessian is represented by
$-DT_h(\psi)$. Its condition number is
\[
\kappa(-DT_h(\psi))
=
\frac{\lambda_{\max}(-DT_h(\psi))}
     {\lambda_{\min}(-DT_h(\psi))}.
\]
It increases
from approximately $92.2$ initially to $863.2$ near convergence,
where the eigenvalues range from $7.02$ to $6063.30$. This spread explains the slow gradient convergence that we observe. Near an optimum,
a Hessian eigenmode with eigenvalue $\lambda$ is multiplied by
$1-\eta_{\mathrm{GF}}\lambda$ at each gradient step.
Local linear stability then requires
\[
0<\eta_{\mathrm{GF}}
<\frac{2}{\lambda_{\max}}
\simeq3.30\times10^{-4},
\]
and the modes associated with small eigenvalues then decay slowly, as we must choose a very small step size $\eta_{\mathrm{GF}}$. For mirror flow, exact inverse-Hessian preconditioning gives
$(-DT_h(\psi))^{-1}(-DT_h(\psi))=I$, whose condition number is $1$.
The \emph{preconditioned} condition number is therefore small,
and the local mirror iteration multiplies every mode by
$1-\eta_{\mathrm{MF}}$. Here, we can choose a larger time step.
By compensating for the Hessian’s uneven curvature, preconditioning therefore permits a larger stable outer optimization step which explains mirror flow’s faster convergence.

\section{Proofs of the mirror-flow and terminal-map results}
\label{sec:proofs-mirror-terminal}

The proofs concerning attainment, the dual value map, and the Bregman
scheme retain the general Markovian cost $L=L(t,x,a,b)$ introduced in
the preliminaries.  The autonomous notation
\eqref{eq:autonomous_H_control_representation} applies from the proof of
Theorem~\ref{prop:frechet_terminal_map} onward.

\subsection{Proofs of Propositions~\ref{prop:attainment_A} and~\ref{prop:A_diff}}

\begin{proof}[Proof of Proposition~\ref{prop:attainment_A}]
The argument is based on the compactness, stability, and
lower-semicontinuity results of Tan--Touzi; see in particular
Lemma~3.9, Proposition~3.11, Remark~3.12, and the proof of
Lemma~3.13 therein.

Set
$
J_\psi(P):=\mathcal C(P)+\langle\psi,P_1\rangle$ and
$
\mathcal C(P):=
\mathbb E^P\left[\int_0^1 L(t,X_t,a_t^P,b_t^P)\,dt\right].
$
Since $\psi\in C(\mathbb T^d)$, $L\geq0$, and there exists at least
one admissible law with finite running cost (see Remark~\ref{rem:finite-cost-feasibility}), we have
$-\|\psi\|_\infty\leq \mathcal A(\psi)<\infty$.

\noindent
Let $(P_n)$ be a minimizing sequence for $\mathcal A(\psi)$. After discarding
finitely many terms, we may assume that
$J_\psi(P_n)\leq \mathcal A(\psi)+1$. Since
$\langle\psi,(P_n)_1\rangle
\geq-\|\psi\|_\infty$, it follows that
$
\sup_n\mathcal C(P_n)
\leq \mathcal A(\psi)+1+\|\psi\|_\infty<\infty.
$
The coercivity assumption therefore yields
$
\sup_n\mathbb E^{P_n}\int_0^1
\left(|a_t^{P_n}|^p+|b_t^{P_n}|^p\right)\,dt<\infty.
$

For each $n$, let $\overline P_n$ be the probability measure on the
enlarged canonical space induced by $(X,A^{P_n},B^{P_n})$, where
$A_t^{P_n}:=\int_0^t a_s^{P_n}\,\mathrm ds$ and
$B_t^{P_n}:=\int_0^t b_s^{P_n}\,\mathrm ds$. Let $\pi:\mathbb R^d\to\mathbb T^d$ be the quotient map. For the
application of~\cite{TaTo13}, lift $X$ continuously to $\widetilde X$,
choosing $\widetilde X_0$ as the representative of $X_0$ in $[0,1)^d$.
The compactness and projection arguments apply to
$(\widetilde X,2A,B)$ with control set $\{(2a,b):(a,b)\in U\}$ and
running cost $(t,y,\alpha,b)\mapsto L(t,\pi(y),\alpha/2,b)$.
The torus path and its lift generate the same filtration under this
choice of initial representative. In our coordinates the enlarged
cost is
$
\overline{\mathcal C}(\overline P)
:=\mathbb E^{\overline P}\!\int_0^1
L(t,X_t,\dot A_t,\dot B_t)\,\mathrm dt,
$
with drift characteristic \(B\) and covariance characteristic \(2A\), where \((\dot A_t,\dot B_t)\in U\) almost everywhere.
Thus $\overline{\mathcal C}(\overline P_n)=\mathcal C(P_n)$.
\noindent
The preceding $L^p$ estimate allows us to apply the tightness argument
used in the proof of Tan--Touzi, Lemma~3.13. Hence, after extraction
of a subsequence, still denoted by $(\overline P_n)$,
$\overline P_n\Longrightarrow\overline P$ on the enlarged canonical
space. The stability argument in the same proof shows that, under
$\overline P$, the canonical process $X$ remains a semimartingale
with absolutely continuous characteristics $(A,B)$ whose density
$(a_t,b_t)$ belongs to $U$ for
$d\overline P\otimes dt$-almost every $(\omega,t)$. Moreover, since
the map $\omega\mapsto X_0(\omega)$ is continuous and $(P_n)_0=\mu$
for every $n$, the initial marginal of $\overline P$ is also $\mu$.
Thus $\overline P$ is admissible in the enlarged formulation.

To pass to the limit in the running cost, write $Z=(A,B)$ and let
$
I(x,Z):=\int_0^1 L(t,x_t,\dot Z_t)\,\mathrm dt
$ when $Z$ is absolutely continuous with $\dot Z_t\in U$ almost
everywhere, and set $I(x,Z)=+\infty$ otherwise. We claim that $I$
is lower semicontinuous for uniform convergence.

Indeed, let $(x_n,Z_n)\to(x,Z)$ uniformly and take a subsequence
realizing a finite $\liminf_n I(x_n,Z_n)$. Coercivity bounds
$\dot Z_n$ in $L^p$, so, after extraction,
$\mathrm dt\,\delta_{\dot Z_n(t)}(\mathrm du)$ converges weakly
to a measure of the form $\mathrm dt\,\lambda_t(\mathrm du)$
on $[0,1]\times U$. Since $p>1$, the first moments are uniformly
integrable. Testing against smooth functions of time and integrating
by parts identifies the weak derivative of $Z$ as
$\dot Z_t=\int_U u\,\lambda_t(\mathrm du)\in U$ almost everywhere.
In particular, $Z$ is absolutely continuous. By uniform convergence
of $x_n$, the joint measures
$\mathrm dt\,\delta_{x_n(t)}(\mathrm dy)
\delta_{\dot Z_n(t)}(\mathrm du)$
converge to $\mathrm dt\,\delta_{x(t)}(\mathrm dy)
\lambda_t(\mathrm du)$. Portmanteau and Jensen's inequality give
\[
\liminf_{n\to\infty} I(x_n,Z_n)
\geq \int_0^1\!\int_U L(t,x_t,u)\,\lambda_t(\mathrm du)\,\mathrm dt
\geq I(x,Z).
\]
Applying Portmanteau on the enlarged canonical space now yields
$
\overline{\mathcal C}(\overline P)
\leq \liminf_{n\to\infty}\overline{\mathcal C}(\overline P_n)
= \liminf_{n\to\infty}\mathcal C(P_n).
$
Since $(X,A,B)\mapsto\psi(X_1)$ is bounded and continuous,
$\mathbb E^{\overline P_n}[\psi(X_1)]
\to\mathbb E^{\overline P}[\psi(X_1)]$. As $\overline P_n$ and $P_n$
have the same $X$-marginal,
$\mathbb E^{\overline P_n}[\psi(X_1)]
=\langle\psi,(P_n)_1\rangle$. Consequently,

$$
\overline{\mathcal C}(\overline P)
+\mathbb E^{\overline P}[\psi(X_1)]
\leq
\liminf_{n\to\infty}
\left(\mathcal C(P_n)+\langle\psi,(P_n)_1\rangle\right)
=\mathcal A(\psi).
$$

In particular, $\overline{\mathcal C}(\overline P)<\infty$. Set
$\nu^\ast:=\overline P\circ X_1^{-1}$, so that
$\overline P\in\overline{\mathcal P}(\mu,\nu^\ast)$. By the
coercivity assumption and Tan--Touzi, Remark~3.12,
$
\mathbb E^{\overline P}
\int_0^1 \bigl(|a_s|+|b_s|\bigr)\,\mathrm ds
<\infty.
$ Hence
Proposition~3.11(ii) yields a law
$P\in\mathcal P(\mu,\nu^\ast)$ such that
$\mathcal C(P)\leq\overline{\mathcal C}(\overline P)$. In particular,
$P_1=\nu^\ast=\overline P\circ X_1^{-1}$, and therefore
$\langle\psi,P_1\rangle
=\mathbb E^{\overline P}[\psi(X_1)]$.
\noindent
Since $P$ is admissible in the definition of $\mathcal A(\psi)$,

$$
\mathcal A(\psi)
\leq\mathcal C(P)+\langle\psi,P_1\rangle
\leq\overline{\mathcal C}(\overline P)
+\mathbb E^{\overline P}[\psi(X_1)]
\leq \mathcal A(\psi).
$$
\noindent
All inequalities are therefore equalities, so
$
\mathcal A(\psi)=\mathcal C(P)+\langle\psi,P_1\rangle.
$
Thus $P$ is a minimizer of $\mathcal A(\psi)$.
\end{proof}

\begin{proof}[Proof of Proposition \ref{lem:terminal_marginal_uniqueness}]
Fix a terminal potential $\psi$, and let $P^1$ and $P^2$ be two
minimizers of $\mathcal A(\psi)$. Denote their terminal marginals by $
\nu_1:=(P^1)_1,
$ and
$\nu_2:=(P^2)_1.
$

We first claim that
$
\mathcal C(P^i)=V(\nu_i),
i=1,2.
$
Indeed, by definition of $V$, $
V(\nu_i)\leq \mathcal C(P^i).
$
Suppose, for instance, that
$
V(\nu_1)<\mathcal C(P^1).
$
Set
$
\varepsilon
:=
\frac{\mathcal C(P^1)-V(\nu_1)}{2}>0.
$
By the definition of the infimum, there exists an admissible law
$Q\in\mathcal P(\mu)$ with $Q_1=\nu_1$ such that
$
\mathcal C(Q)
<
V(\nu_1)+\varepsilon
=
\frac{V(\nu_1)+\mathcal C(P^1)}{2}
<
\mathcal C(P^1).
$
\noindent
Since $Q_1=(P^1)_1=\nu_1$, we also have
$\langle\psi,Q_1\rangle=\langle\psi,\nu_1\rangle$. Hence
$
\mathcal C(Q)+\langle\psi,Q_1\rangle
<
\mathcal C(P^1)+\langle\psi,(P^1)_1\rangle
=
\mathcal A(\psi),
$
contradicting the optimality of $P^1$. Therefore
$\mathcal C(P^1)=V(\nu_1)$, and the same argument gives
$\mathcal C(P^2)=V(\nu_2)$.
Consequently,
$
\mathcal A(\psi)
=
V(\nu_i)+\langle\psi,\nu_i\rangle,
i=1,2.
$

Suppose now, by contradiction, that $\nu_1\neq\nu_2$. Since $V$ is
strictly convex on its effective domain,
$
V\left(\frac{\nu_1+\nu_2}{2}\right)
<
\frac12V(\nu_1)+\frac12V(\nu_2).
$
\noindent
By linearity of the pairing,
$
\left\langle
\psi,\frac{\nu_1+\nu_2}{2}
\right\rangle
=
\frac12\langle\psi,\nu_1\rangle
+
\frac12\langle\psi,\nu_2\rangle.
$
\noindent
Adding the two relations yields
$
V\left(\frac{\nu_1+\nu_2}{2}\right)
+
\left\langle
\psi,\frac{\nu_1+\nu_2}{2}
\right\rangle
<
\frac12
\left(
V(\nu_1)+\langle\psi,\nu_1\rangle
\right)
+
\frac12
\left(
V(\nu_2)+\langle\psi,\nu_2\rangle
\right)
=
\mathcal A(\psi).
$
\noindent
But, by definition,
$
\mathcal A(\psi)
=
\inf_{\nu}
\left\{
V(\nu)+\langle\psi,\nu\rangle
\right\},
$
which is a contradiction.
Therefore $\nu_1=\nu_2$. Hence all minimizers of $\mathcal A(\psi)$ have the
same terminal marginal, which we denote by $T(\psi)$.
\end{proof}

\begin{lemma}[Local compactness and stability of optimal laws]
\label{lem:local_compactness_optimizers}

Let Assumption~\ref{ass:attainment_uniqueness}
hold. Fix $\psi,h\in C(\mathbb T^d)$
and $\delta>0$. For any sequence
$\varepsilon_n\in[-\delta,\delta]$
and any choice of optimal laws
$P_n$ for $\mathcal A(\psi+\varepsilon_nh)$,
the sequence $(P_n)$ is relatively
compact for weak convergence.
\noindent
More precisely, after extraction of
a subsequence, the corresponding laws on the enlarged canonical
space converge weakly to a law which can be projected onto an
admissible law $P$ on the original canonical space, and
$
\mathcal C(P)
\leq
\liminf_{n\to\infty}\mathcal C(P_n).
$

\end{lemma}

\begin{proof}
Let $P^0$ be an admissible law with $\mathcal C(P^0)<\infty$,
whose existence follows from Remark~\ref{rem:finite-cost-feasibility}.
Since $P_n$ minimizes $\mathcal A(\psi+\varepsilon_n h)$, optimality gives
$
\mathcal C(P_n)
+\langle\psi+\varepsilon_n h,(P_n)_1\rangle
\leq
\mathcal C(P^0)
+\langle\psi+\varepsilon_n h,(P^0)_1\rangle.
$
As both terminal marginals are probability measures and
$|\varepsilon_n|\leq\delta$, we obtain
\[
\mathcal C(P_n)
\leq \mathcal C(P^0)+2\|\psi+\varepsilon_n h\|_\infty
\leq \mathcal C(P^0)+2\bigl(\|\psi\|_\infty+\delta\|h\|_\infty\bigr).
\]
By coercivity,
$
\sup_n\mathbb E^{P_n}\int_0^1
\bigl(|a_t^{P_n}|^p+|b_t^{P_n}|^p\bigr)\,\mathrm dt<\infty.
$

Let $\overline P_n$ be the law of $(X,A^{P_n},B^{P_n})$ under $P_n$,
so that $\overline{\mathcal C}(\overline P_n)=\mathcal C(P_n)$.
The uniform $L^p$ bound gives tightness of $(\overline P_n)$ by the
same lift and compactness argument as in the proof of
Proposition~\ref{prop:attainment_A}. Since projection onto the
$X$-coordinate is continuous, $(P_n)$ is relatively compact.

Choose a subsequence along which $\mathcal C(P_n)$ converges to its
liminf, and extract a further subsequence $(n_k)$ such that
$\overline P_{n_k}\Longrightarrow\overline P$. The stability argument
from the same proof gives absolutely continuous characteristics
with $(\dot A_t,\dot B_t)\in U$ for
$\mathrm d\overline P\otimes\mathrm dt$-almost every $(\omega,t)$.
Continuity of the initial evaluation map gives
$\overline P\circ X_0^{-1}=\mu$. The lower semicontinuity of the
enlarged cost established in that proof yields
\[
\overline{\mathcal C}(\overline P)
\leq \liminf_{k\to\infty}\mathcal C(P_{n_k})
= \liminf_{n\to\infty}\mathcal C(P_n)<\infty.
\]
By coercivity and the projection argument in the proof of
Proposition~\ref{prop:attainment_A}, the law
$P=\overline P\circ X^{-1}$ belongs to $\mathcal P(\mu)$ and satisfies
$
\mathcal C(P)\leq\overline{\mathcal C}(\overline P)
\leq\liminf_{n\to\infty}\mathcal C(P_n).
$
Moreover, $P_{n_k}\Longrightarrow P$ by continuity of projection onto
the $X$-coordinate.
\end{proof}
\begin{lemma}[Weak continuity of the optimal terminal marginal]
\label{lem:weak_continuity_T}
Let Assumption~\ref{ass:attainment_uniqueness} hold. Then, for every
$\psi,h\in C(\mathbb T^d)$,
$
T(\psi+\varepsilon h)\Longrightarrow T(\psi)$
$\text{as }\varepsilon\to0.
$
\end{lemma}

\begin{proof}
Let $\varepsilon_n\to0$, and for each $n$ choose an optimizer $P_n$
of $\mathcal A(\psi+\varepsilon_n h)$. By
Lemma~\ref{lem:local_compactness_optimizers}, every subsequence of
$(P_n)$ admits a further subsequence, still denoted by $(P_n)$, which
converges weakly to an admissible law $P$ satisfying
$
\mathcal C(P)
\leq
\liminf_{n\to\infty}\mathcal C(P_n).
$
\noindent
Let $Q$ be an optimizer of $\mathcal A(\psi)$. By optimality of $P_n$,

$$
\mathcal C(P_n)
+\langle\psi+\varepsilon_nh,(P_n)_1\rangle
\leq
\mathcal C(Q)
+\langle\psi+\varepsilon_nh,Q_1\rangle.
$$
\noindent
Since $P_n\Longrightarrow P$, the continuity of the terminal
evaluation map gives $(P_n)_1\Longrightarrow P_1$. As
$\psi,h\in C(\mathbb T^d)$,

$$
\langle\psi,(P_n)_1\rangle\to\langle\psi,P_1\rangle,
\qquad
\varepsilon_n\langle h,(P_n)_1\rangle\to0.
$$
\noindent
Passing to the limit inferior in the optimality inequality therefore
yields
$
\mathcal C(P)+\langle\psi,P_1\rangle
\leq
\mathcal C(Q)+\langle\psi,Q_1\rangle
=
\mathcal A(\psi).
$
\noindent
Since $P$ is admissible, the reverse inequality follows from the
definition of $\mathcal A(\psi)$. Hence $P$ is an optimizer of $\mathcal A(\psi)$.
By Proposition~\ref{lem:terminal_marginal_uniqueness}, $
P_1=T(\psi).
$
Moreover, since $P_n$ is optimal for $\mathcal A(\psi+\varepsilon_n h)$, we have $
(P_n)_1=T(\psi+\varepsilon_n h).
$
Thus every subsequence of
$\bigl(T(\psi+\varepsilon_n h)\bigr)_n$ admits a further subsequence
converging weakly to $T(\psi)$. Hence the whole sequence converges:
$
T(\psi+\varepsilon_n h)\Longrightarrow T(\psi).
$
\noindent
Since this holds for every sequence $\varepsilon_n\to0$, we conclude
that
$
T(\psi+\varepsilon h)\Longrightarrow T(\psi)$
$\text{as }\varepsilon\to0.$
\end{proof}

\begin{proof}[Proof of Proposition \ref{prop:A_diff}]
We use the abstract stochastic-control representation \eqref{eq:stoch_A}.
Let $P^\psi$ be optimal for $\psi$. Using $P^\psi$ as a competitor for the problem with terminal payoff $\psi+\varepsilon h$ gives
\[
\begin{aligned}
\mathcal A(\psi+\varepsilon h)
&\le
\mathcal C(P^\psi)
+
\pair{\psi+\varepsilon h}{P^\psi_1}
=
\mathcal A(\psi)
+
\varepsilon\pair{h}{T(\psi)}.
\end{aligned}
\]
Therefore
$
\mathcal A(\psi+\varepsilon h)-\mathcal A(\psi)
\le
\varepsilon\pair{h}{T(\psi)}.
$
Conversely, let $P^{\psi+\varepsilon h}$ be optimal for $\psi+\varepsilon h$.
Then
\[
\begin{aligned}
\mathcal A(\psi+\varepsilon h)
&=
\mathcal C(P^{\psi+\varepsilon h})
+
\pair{\psi+\varepsilon h}{P^{\psi+\varepsilon h}_1}      
=
\mathcal C(P^{\psi+\varepsilon h})
+
\pair{\psi}{P^{\psi+\varepsilon h}_1}
+
\varepsilon\pair{h}{P^{\psi+\varepsilon h}_1}.
\end{aligned}
\]
Since $P^{\psi+\varepsilon h}$ is admissible for the problem with terminal payoff $\psi$,
$
\mathcal C(P^{\psi+\varepsilon h})
+
\pair{\psi}{P^{\psi+\varepsilon h}_1}
\ge
\mathcal A(\psi).
$
Hence
$
\mathcal A(\psi+\varepsilon h)-\mathcal A(\psi)
\ge
\varepsilon\pair{h}{T(\psi+\varepsilon h)}.
$
For $\varepsilon>0$, we obtain
$$
\pair{h}{T(\psi+\varepsilon h)}
\le
\frac{\mathcal A(\psi+\varepsilon h)-\mathcal A(\psi)}{\varepsilon}
\le
\pair{h}{T(\psi)}.
$$
By Lemma~\ref{lem:weak_continuity_T}, letting $\varepsilon\downarrow0$ gives
$
\lim_{\varepsilon\downarrow0}
\frac{\mathcal A(\psi+\varepsilon h)-\mathcal A(\psi)}{\varepsilon}
=
\pair{h}{T(\psi)}.
$
The same argument with $\varepsilon<0$ gives the left derivative. Therefore
$
D\mathcal A(\psi)[h]
=
\pair{T(\psi)}{h}.
$
The formula for $D\Phi$ follows from $\Phi=-\mathcal A$.
\end{proof}

\subsection{Proof of Theorem~\ref{thm:Bregman_contraction}}
\begin{proof}
For any $\bar \psi_0,\bar\psi_1,\bar\psi_2\in\calX$, we have 
$$
\calB_\Phi(\bar\psi_0\mid\bar\psi_2)
=
\calB_\Phi(\bar\psi_0\mid
\bar\psi_1)
+
\calB_\Phi(\bar\psi_1\mid\bar\psi_2)
+
\pair{
D\Phi(\bar\psi_1)-D\Phi(\bar\psi_2)
}{
\bar\psi_0-\bar\psi_1
}.
$$
Apply it with
$\bar\psi_0=\psi^\star$,
$\bar\psi_1=\bar\psi_{k+1}$,
$\bar\psi_2=\bar\psi_k$.
We get
\[
\begin{aligned}
\calB_\Phi(\psi^\star\mid\bar\psi_k)
&=
\calB_\Phi(\psi^\star\mid\bar\psi_{k+1})
+
\calB_\Phi(\bar\psi_{k+1}\mid\bar\psi_k)     
+
\pair{
D\Phi(\bar\psi_{k+1})-D\Phi(\bar\psi_k)
}{
\psi^\star-\bar\psi_{k+1}
}.
\end{aligned}
\]
The discrete mirror descent equation gives
$
D\Phi(\bar\psi_{k+1})-D\Phi(\bar\psi_k)
=
-\tau
\left(
D\Phi(\bar\psi_{k+1})-D\Phi(\psi^\star)
\right).
$
Therefore
\[
\begin{aligned}
\pair{
D\Phi(\bar\psi_{k+1})-D\Phi(\bar\psi_k)
}{
\psi^\star-\bar\psi_{k+1}
}                                                
&=
\tau
\pair{
D\Phi(\bar\psi_{k+1})-D\Phi(\psi^\star)
}{
\bar\psi_{k+1}-\psi^\star
}                                          \\
&= 
\tau
\left[
\calB_\Phi(\bar\psi_{k+1}\mid\psi^\star)
+
\calB_\Phi(\psi^\star\mid\bar\psi_{k+1})
\right].
\end{aligned}
\]
Substituting this into the previous identity gives
\[
\begin{aligned}
\calB_\Phi(\psi^\star\mid\bar\psi_k)
&=
(1+\tau)\calB_\Phi(\psi^\star\mid\bar\psi_{k+1})
+
\calB_\Phi(\bar\psi_{k+1}\mid\bar\psi_k)           
+
\tau\calB_\Phi(\bar\psi_{k+1}\mid\psi^\star).
\end{aligned}
\]
Since all Bregman divergences are nonnegative,
$
(1+\tau)\calB_\Phi(\psi^\star\mid\bar\psi_{k+1})
\le
\calB_\Phi(\psi^\star\mid\bar\psi_k).
$
Iterating yields the result.
\end{proof}

\subsection{Proof of Corollary~\ref{cor:value_convergence}}

\begin{proof}
Since $T(\psi^\star)=\nu$, we have
$
D\Phi(\psi^\star)=-\nu.
$
Hence
\begin{align*}
\calB_\Phi(\psi\mid\psi^\star)
&=
\Phi(\psi)-\Phi(\psi^\star)
-
D\Phi(\psi^\star)[\psi-\psi^\star] \\
&=
\Phi(\psi)-\Phi(\psi^\star)
+
\pair{\nu}{\psi-\psi^\star} \\
&=
D_\nu(\psi^\star)-D_\nu(\psi).
\end{align*}
On the other hand, the discrete Bregman identity gives
$
\tau
\calB_\Phi(\bar\psi_{k+1}\mid\psi^\star)
\leq
\calB_\Phi(\psi^\star\mid\bar\psi_k).
$
Combining this with the reverse-Bregman contraction,
$
\calB_\Phi(\psi^\star\mid\bar\psi_k)
\leq
(1+\tau)^{-k}
\calB_\Phi(\psi^\star\mid\bar\psi_0),
$
yields
\[
D_\nu(\psi^\star)-D_\nu(\bar\psi_{k+1})
=
\calB_\Phi(\bar\psi_{k+1}\mid\psi^\star)
\leq
\frac{1}{\tau}(1+\tau)^{-k}
\calB_\Phi(\psi^\star\mid\bar\psi_0),
\]
which proves \eqref{eq:value_convergence}.
\end{proof}

\subsection{Proof of Theorem~\ref{prop:frechet_terminal_map}}

We establish the differentiability of the
terminal map in three steps. First, we prove
well-posedness of the nonlinear HJB--Fokker--Planck
system. Next, we differentiate the HJB equation
and the optimal feedback. Finally, we show
that the linearized Fokker--Planck equation
has a unique solution and identify the
derivative of the terminal map.

For
$\psi,h\in C^{m+2,\alpha}_\diamond(\mathbb T^d)$,
the linearized system is
\begin{equation}
\left\{
\begin{aligned}
& \partial_t w^h
+b^\psi\cdot\nabla w^h
+a^\psi:\nabla^2 w^h =0, \quad
 &&w^h(1,\cdot)=h,\\
& D^2 L(a^\psi,b^\psi)
\big[(\dot a^h,\dot b^h),(A,\beta)\big]
+\nabla^2 w^h:A
+\nabla w^h\cdot\beta =0,
 &&\forall\,(A,\beta)\in\mathbb E^d,\\
& \partial_t q^h
=-\partial_i\!\left(b_i^\psi q^h+\dot b_i^h m^\psi\right)
+\partial_{ij}\!\left(a_{ij}^\psi q^h+\dot a_{ij}^h m^\psi\right),
&& q^h(0,\cdot)=0.
\end{aligned}
\right.
\label{eq:linearized_system}
\end{equation}
We start by discussing the existence and regularity of the coefficients of this system.
\begin{lemma}[Well-posedness and regularity of the HJB–Fokker–Planck system]
\label{lem:state_system_terminal_map}
\label{rem:reg_HJB}
Under Assumption~\ref{ass:T_diff}, the following assertions hold for
every $\psi\in C^{m+2,\alpha}_\diamond(\mathbb T^d)$.
\begin{enumerate}
\item
Equation~\eqref{eq:HJB_T_differentiability} has a unique solution
$u^\psi\in C^{m+2,\alpha}(Q)$.

\item
The minimizer map
$\mathfrak z=(\mathfrak a,\mathfrak b)$
defined by
\eqref{eq:autonomous_H_control_representation}
belongs to
$C^{m+2,\alpha}_{\mathrm{loc}}
(\mathbb E^d;\mathbb E^d)$.
Moreover, the envelope identities
\begin{equation}
\label{eq:envelope_identity}
D_MH(M,p)=\mathfrak a(M,p),
\qquad
D_pH(M,p)=\mathfrak b(M,p)
\end{equation}
hold for every $(M,p)\in\mathbb E^d$.
Consequently,
\[
(a^\psi,b^\psi)\in
C^{m,\alpha}(Q;\mathbb E^d),
\qquad
\lambda I_d\leq a^\psi\leq\Lambda I_d
\quad\hbox{on }Q,
\]
where $\lambda$ and $\Lambda$ are the
ellipticity constants from
Assumption~\ref{ass:T_diff}\textup{(iv)}.

\item
Equation~\eqref{eq:FP_T_differentiability} has a unique solution
$m^\psi\in C^{m,\alpha}(Q)$.  Moreover,
\[
m^\psi\geq0,
\qquad
\int_{\mathbb T^d}m^\psi(t,x)\,\mathrm dx=1,
\qquad 0\leq t\leq1,
\]
and
\begin{equation}
\label{eq:strict_positivity_density}
m^\psi(t,x)>0,
\qquad (t,x)\in(0,1]\times\mathbb T^d.
\end{equation}
\end{enumerate}
\end{lemma}

\begin{proof}
$ $\newline
\noindent\emph{Minimizer map and Hamiltonian (part of item 2).}
We first establish the regularity of the optimizer
map and of the Hamiltonian.

Fix an arbitrary open set
$\mathcal P\Subset\mathbb E^d$.
By Assumption~\ref{ass:T_diff}\textup{(ii)},
the optimizer $\mathfrak z(M,p)$ is unique and
takes values in a compact subset
$\mathcal K_{\mathcal P}\Subset U$ for
$(M,p)\in\mathcal P$.

The optimizer map is continuous. Indeed, if
$(M_n,p_n)\to(M,p)$ in $\mathcal P$, the
corresponding optimizers belong to the compact
set $\mathcal K_{\mathcal P}$. Every cluster
point minimizes the limiting objective and
therefore equals the unique optimizer
$\mathfrak z(M,p)$. Thus the whole sequence
converges.

Since $\mathfrak z(M,p)\in\operatorname{int}U$,
the first-order optimality condition gives
\begin{equation}
\label{eq:optimality_condition}
DL\bigl(\mathfrak z(M,p)\bigr)+(M,p)=0.
\end{equation}
By Assumption~\ref{ass:T_diff}\textup{(iii)},
$D^2L(\mathfrak z(M,p))$ is invertible.
The finite-dimensional implicit function
theorem therefore provides a local
$C^{m+2,\alpha}$ branch of the optimizer.
By continuity and uniqueness, these local
branches agree with $\mathfrak z$ and hence
\[
\mathfrak z\in
C^{m+2,\alpha}_{\mathrm{loc}}
(\mathbb E^d;\operatorname{int}U).
\]

Differentiating
$
H(M,p)
=
\mathfrak z(M,p)\boldsymbol\cdot(M,p)
+L\bigl(\mathfrak z(M,p)\bigr)
$
in the direction $(N,q)$ and using
\eqref{eq:optimality_condition}, we obtain
\[
\begin{aligned}
DH(M,p)[N,q]
&=
\mathfrak z(M,p)\boldsymbol\cdot(N,q)
\\
&\quad+
\bigl[
DL(\mathfrak z(M,p))+(M,p)
\bigr]\boldsymbol\cdot
D\mathfrak z(M,p)[N,q]
\\
&=
\mathfrak a(M,p):N+
\mathfrak b(M,p)\cdot q.
\end{aligned}
\]
This proves \eqref{eq:envelope_identity}.
In particular,
$
DH=\mathfrak z$,
$H\in C^{m+3,\alpha}_{\mathrm{loc}}(\mathbb E^d)$.


\medskip
\noindent\emph{HJB well-posedness (item 1).}
We can now solve the HJB equation.
Set $v(s,x):=u^\psi(1-s,x)$. Then
\begin{equation}
\label{eq:forward_fully_nonlinear_state}
\partial_sv=H(\nabla^2v,\nabla v),
\qquad
v(0,\cdot)=\psi.
\end{equation}
This is a concave uniformly parabolic Cauchy problem, with a
gradient-dependent lower-order term, on the closed manifold
$\mathbb T^d$. The local \(C^{1+\alpha/2,\,2+\alpha}\) existence and
uniqueness theorem \cite[Theorem~1.2]{Huang2015Cauchy} applies because
$H$ is $C^2$ in a neighborhood of the
compact set
$\bigl\{(\nabla^2\psi(x),\nabla\psi(x)):
x\in\mathbb T^d\bigr\}$, and the eigenvalues of the
principal coefficient $D_MH$ lie in $[\lambda,\Lambda]$.  Since
$\psi\in C^{m+2,\alpha}$ and $H\in C^{m+3,\alpha}_{\mathrm{loc}}$,
successive spatial differentiation and the periodic Schauder estimate
\cite[Chapter~IV, Theorem~5.1]{LadyzhenskayaSolonnikovUraltseva1968}
upgrade the local solution to the regularity claimed in part~1.

For continuation to $s=1$, set $z:=\partial_sv$ and
\[
 \widehat a:=D_MH(\nabla^2v,\nabla v),
 \qquad
 \widehat b:=D_pH(\nabla^2v,\nabla v).
\]
For a constant unit vector $e\in\mathbb R^d$,
write
$ v_e:=e\cdot\nabla v $, $v_{ee}:=e^\top\nabla^2v\,e$.
Differentiating
\eqref{eq:forward_fully_nonlinear_state} gives
$
\partial_sz-\widehat b\cdot\nabla z
-\widehat a:\nabla^2z=0.
$
The maximum principle therefore yields
\begin{equation}
\label{eq:state-time-derivative-bound}
\|\partial_sv(s)\|_{L^\infty}
\leq
\|H(\nabla^2\psi,\nabla\psi)\|_{L^\infty}
\end{equation}
on the existence interval.  For every constant unit vector $e$,
$
 (\partial_s-\widehat b\cdot\nabla
 -\widehat a:\nabla^2)v_e=0.
$
Applying the maximum principle to $v_e$
and $-v_e$, and then taking the supremum
over all unit vectors $e$, gives
$
\|\nabla v(s,\cdot)\|_{L^\infty(\mathbb T^d)}
\leq
\|\nabla\psi\|_{L^\infty(\mathbb T^d)}
$
on the existence interval. A second differentiation
and joint concavity of $H$ in $(M,p)$ give
\[
\bigl(\partial_s-\widehat b\cdot\nabla
-\widehat a:\nabla^2\bigr)v_{ee}
=D^2H(\nabla^2v,\nabla v)
[\J v_e,\J v_e]
\leq0.
\]
The maximum principle and the supremum over $e$ therefore bound the
largest eigenvalue of $\nabla^2v$ by
\[
K_+:=\max\left\{0,
\sup_{x\in\mathbb T^d}\lambda_{\max}(\nabla^2\psi(x))\right\}.
\]
For the negative eigenvalues, write $M=M_+-M_-$ with $M_\pm\geq0$.
Uniform ellipticity gives
\[
H(M,p)-H(0,p)
=\int_0^1D_MH(\theta M,p):M\,\mathrm d\theta
\leq \Lambda\operatorname{Tr}M_+
-\lambda\operatorname{Tr}M_-.
\]
With $p=\nabla v(s,x)$ and $M=\nabla^2v(s,x)$, the gradient bound
controls $H(0,\nabla v)$, \eqref{eq:state-time-derivative-bound}
controls the left-hand side from below, and
$\operatorname{Tr}M_+\leq dK_+$.  Hence
$\operatorname{Tr}M_-$, and thus $\nabla^2v$ from both sides, is
uniformly bounded.

Let $S\in(0,\infty]$ denote the maximal existence time of the forward
solution $v$ of \eqref{eq:forward_fully_nonlinear_state}.
Fix $0<\varepsilon<S$. For every $S'\in(\varepsilon,S)$, the
gradient-dependent parabolic Evans--Krylov estimate
\cite[Theorem~4]{Andrews2004FullyNonlinear}
applies on $[\varepsilon,S']$ in finitely many periodic charts.
The preceding bounds make the estimate uniform in $S'$. Indeed, $H$ is concave in the Hessian
variable and uniformly elliptic.
The bounds established above on
$\nabla v$ and $\nabla^2v$ keep
$\mathcal Jv$ in a fixed compact subset
of $\mathbb E^d$, on which the derivatives
of $H$ are bounded. The preceding
estimates and periodicity provide the required $C^2$ bounds.  We obtain
a uniform $C^{2,\gamma}$ estimate for some
$\gamma\in(0,\alpha]$, after lowering the exponent if necessary.
Schauder estimates for the differentiated equations then give uniform
$C^{m+3,\gamma}$ bounds away from $s=0$.  Since
$m+3+\gamma>m+2+\alpha$, these imply the required
$C^{m+2,\alpha}$ bounds.

If the maximal existence time $S$ were finite, the equation would also
bound $\partial_s v=H(\nabla^2v,\nabla v)$ uniformly in
$C^{m+1,\gamma}$ on $[\varepsilon,S)$.  Thus $v(s)$ would be Cauchy in
$C^{m+1,\gamma}$ as $s\uparrow S$, and interpolation with the uniform
$C^{m+3,\gamma}$ bound would make it Cauchy in $C^{m+2,\alpha}$ because
$m+2+\alpha<m+3+\gamma$.  Its $C^{m+2,\alpha}$ trace at $S$ would allow
Theorem~1.2 of \cite{Huang2015Cauchy} to restart the solution, a
contradiction.  Hence the solution exists on $[0,1]$, and the same
differentiated Schauder estimates give
$v\in C^{1+(m+\alpha)/2,m+2+\alpha}(Q)$.  

Uniqueness follows from the comparison
principle. Indeed, if $v_1$ and $v_2$
are two solutions with the same initial
datum, their difference $w:=v_1-v_2$
satisfies
\[
\partial_s w
=\widetilde a:\nabla^2w
+\widetilde b\cdot\nabla w,
\qquad w(0,\cdot)=0,
\]
where
\[
\widetilde a
:=
\int_0^1
D_MH\bigl(
\mathcal Jv_2+\theta\mathcal Jw
\bigr)\,\mathrm d\theta, \qquad
\widetilde b
:=
\int_0^1
D_pH\bigl(
\mathcal Jv_2+\theta\mathcal Jw
\bigr)\,\mathrm d\theta.
\]
By Assumption~\ref{ass:T_diff}\textup{(iv)},
$\widetilde a$ is uniformly elliptic.
The maximum principle therefore gives
$w=0$.

\medskip
\noindent\emph{Optimal feedback (completion of item~2).}
By the regularity of $\mathfrak z$
established above and the regularity of
$u^\psi$ proved in item~1, we have
\[
(a^\psi,b^\psi)
=\mathfrak z(\mathcal Ju^\psi)
\in C^{m,\alpha}(Q;\mathbb E^d).
\]
Moreover, the envelope identity
\eqref{eq:envelope_identity} and
Assumption~\ref{ass:T_diff}\textup{(iv)}
give
\[
\lambda I_d
\leq a^\psi
=D_MH(\mathcal Ju^\psi)
\leq\Lambda I_d
\qquad\text{on }Q.
\]

\medskip
\noindent\emph{Fokker–Planck equation (item 3).}
We now turn to the Fokker--Planck equation. Expanding the right-hand side of
\eqref{eq:FP_T_differentiability} gives a linear, uniformly
parabolic equation
\[
\partial_tm^\psi
=a^\psi_{ij}\partial_{ij}m^\psi
+\bigl(2\partial_i a^\psi_{ij}-b^\psi_j\bigr)
 \partial_jm^\psi
+\bigl(\partial_{ij}a^\psi_{ij}
-\partial_i b^\psi_i\bigr)m^\psi.
\]
Set $\ell:=m-2+\alpha>0$.  The principal, first-order, and zeroth-order
coefficients have parabolic H\"older orders at least $\ell$, $\ell+1$,
and $\ell$, respectively, while $\mu\in C^{\ell+2}(\mathbb T^d)$.
After periodic extension to $\mathbb R^d$, the Cauchy-problem Schauder
theorem
\cite[Chapter~IV, Theorem~5.1]{LadyzhenskayaSolonnikovUraltseva1968}
gives a unique solution of parabolic order $\ell+2=m+\alpha$, namely
$m^\psi\in C^{m,\alpha}(Q)$; uniqueness makes it periodic.

Integration over the torus gives mass conservation.  With
$c:=\partial_{ij}a^\psi_{ij}-\partial_i b^\psi_i$ and
$K\geq\|c\|_{L^\infty(Q)}$, the weak maximum principle applied to
$e^{-Kt}m^\psi$, whose zeroth-order coefficient is nonpositive, gives
$m^\psi\geq0$.  Since its initial datum has mass one, the strong maximum
principle (equivalently, the Harnack inequality) in
\cite[Chapter~III, Sections~7 and~10]{LadyzhenskayaSolonnikovUraltseva1968}
yields
\eqref{eq:strict_positivity_density}.
\end{proof}

\begin{remark}[Localization of the optimal feedback]
\label{rem:hessian_localization}
\label{rem:uniform_ellipticity}
For every fixed
$\psi\in C^{m+2,\alpha}_\diamond(\mathbb T^d)$,
the continuity of $\mathcal J u^\psi$ and the
compactness of $Q$ imply that
$\mathcal J u^\psi(Q)$ is compact in
$\mathbb E^d$.
We may therefore choose a bounded open set
$\mathcal P_\psi\Subset\mathbb E^d$
containing this range.
By Assumption~\ref{ass:T_diff}\textup{(ii)},
the associated optimal controls belong to
a compact subset of $\operatorname{int}U$.
Moreover, uniform ellipticity of $a^\psi$
follows from the identity
$
a^\psi=D_MH(\mathcal J u^\psi)
$
and Assumption~\ref{ass:T_diff}\textup{(iv)}.
\end{remark}


\begin{lemma}[Differentiability of the HJB--control map]
\label{lem:HJB_control_differentiability}
Under Assumption~\ref{ass:T_diff}, the map
$\psi\mapsto u^\psi$ is locally of class $C^2$ from
$C^{m+2,\alpha}_\diamond(\mathbb T^d)$ into $C^{m+2,\alpha}(Q)$,
and the map $\psi\mapsto(a^\psi,b^\psi)$ is locally of class
$C^{1,1}$ from the same domain into $C^{m,\alpha}(Q;\mathbb E^d)$.
Their derivatives are
\[
Du^\psi[h]=w^h,
\qquad
D(a^\psi,b^\psi)[h]=\dot z^h,
\]
where
$
w^h\in C^{m+2,\alpha}(Q) $,
$\dot z^h=(\dot a^h,\dot b^h)
\in C^{m,\alpha}(Q;\mathbb E^d)
$
are uniquely determined by the first
two equations of
\eqref{eq:linearized_system}.
\end{lemma}

\begin{proof}

Fix
$\psi\in C^{m+2,\alpha}_\diamond(\mathbb T^d)$.
By Lemma~\ref{lem:state_system_terminal_map},
the solution $u^\psi$ belongs to
$C^{m+2,\alpha}(Q)$.

Choose an open set
$\mathcal P\Subset\mathbb E^d$ containing
$\mathcal J u^\psi(Q)$, and define
$
\mathcal O_\psi
:=
\bigl\{
v\in C^{m+2,\alpha}(Q):
\mathcal Jv(Q)\subset\mathcal P
\bigr\}.
$
This is an open neighborhood of $u^\psi$
in $C^{m+2,\alpha}(Q)$.
Consider
$\mathcal F(v,\phi)
:=
\bigl(
\partial_t v+H(\mathcal Jv),\,
v(1,\cdot)-\phi
\bigr)
$
from
$\mathcal O_\psi
\times C^{m+2,\alpha}_\diamond(\mathbb T^d)$
into
$C^{m,\alpha}(Q)
\times C^{m+2,\alpha}(\mathbb T^d)$.

The nonlinear map
$\mathcal N(v):=H(\mathcal Jv)$
has derivatives
\[
\begin{aligned}
D\mathcal N(v)[z]
&=
D_MH(\mathcal Jv):\nabla^2z
+D_pH(\mathcal Jv)\cdot\nabla z,\\
D^2\mathcal N(v)[z_1,z_2]
&=
D^2H(\mathcal Jv)
[\mathcal Jz_1,\mathcal Jz_2].
\end{aligned}
\]
The local $C^{m+3,\alpha}$ regularity of $H$,
together with the standard H\"older
composition and product estimates,
shows that $\mathcal F$ is of class $C^2$.

By the envelope identities
\eqref{eq:envelope_identity},
\[
D_v\mathcal F(u^\psi,\psi)[z]
=
\bigl(
\partial_tz+b^\psi\cdot\nabla z
+a^\psi:\nabla^2z,\,
z(1,\cdot)
\bigr).
\]
The coefficients $(a^\psi,b^\psi)$
belong to $C^{m,\alpha}(Q;\mathbb E^d)$,
and $a^\psi$ is uniformly elliptic.
The backward Schauder isomorphism theorem
\cite[Chapter~IV, Theorem~5.1]
{LadyzhenskayaSolonnikovUraltseva1968},
after time reversal, therefore gives
\[
D_v\mathcal F(u^\psi,\psi):
C^{m+2,\alpha}(Q)
\xrightarrow{\ \sim\ }
C^{m,\alpha}(Q)
\times C^{m+2,\alpha}(\mathbb T^d).
\]

The Banach-space implicit function theorem,
together with uniqueness of the HJB solution,
shows that $\phi\mapsto u^\phi$ is locally
of class $C^2$ near $\psi$.
For
$h\in C^{m+2,\alpha}_\diamond(\mathbb T^d)$,
differentiating
$\mathcal F(u^\phi,\phi)=0$
at $\phi=\psi$ gives
$w^h:=Du^\psi[h]$, which satisfies
the backward equation in
\eqref{eq:linearized_system}.
The Schauder isomorphism gives uniqueness
of $w^h$ in $C^{m+2,\alpha}(Q)$.

After shrinking the neighborhood of $\psi$
if necessary,
$\mathcal Ju^\phi(Q)\subset\mathcal P$
for every $\phi$ in that neighborhood. Since $\mathfrak z\in C^{m+2,\alpha}_{\mathrm{loc}}$,
composition by $\mathfrak z$ is locally $C^{1,1}$ on
$C^{m,\alpha}(Q;\mathbb E^d)$. Indeed, for jets $Y_1,Y_2$ in a
sufficiently small convex neighborhood of $\mathcal Ju^\psi$,
the H\"older composition and product estimates, applied to the
pointwise mean-value formula, give
\[
\|D\mathfrak z(Y_1)-D\mathfrak z(Y_2)\|_{C^{m,\alpha}(Q)}
\leq C\|Y_1-Y_2\|_{C^{m,\alpha}(Q)}.
\]
The same estimate gives a quadratic remainder for the first-order
expansion of $Y\mapsto\mathfrak z(Y)$, whose derivative is
$K\mapsto D\mathfrak z(Y)[K]$. Since
$\phi\mapsto\mathcal Ju^\phi$ is locally $C^2$, the feedback map
$\phi\mapsto\mathfrak z(\mathcal Ju^\phi)$ is locally $C^{1,1}$.

Set
$\dot z^h:=D(a^\psi,b^\psi)[h]$.
Differentiating the optimality condition
\eqref{eq:optimality_condition} yields
\[
D^2L(a^\psi,b^\psi)
[\dot z^h,(A,\beta)]
+\mathcal Jw^h\boldsymbol\cdot(A,\beta)
=0,
\qquad (A,\beta)\in\mathbb E^d.
\]
Since
$\mathcal Ju^\psi(Q)\subset\mathcal P$,
Assumption~\ref{ass:T_diff}\textup{(ii)--(iii)}
implies that
$D^2L(a^\psi,b^\psi)$ is uniformly
positive definite on $Q$ and has a
uniformly bounded inverse.
 Consequently,
\[
\dot z^h
=
-\bigl(D^2L(a^\psi,b^\psi)\bigr)^{-1}
[\mathcal Jw^h]
\in C^{m,\alpha}(Q;\mathbb E^d).
\]
This is the unique solution of the
algebraic equation in
\eqref{eq:linearized_system}.

We have therefore established the asserted regularity of the HJB and feedback maps and the existence
and uniqueness of $(w^h,\dot z^h)$ in the
stated H\"older spaces.
\end{proof}

\begin{lemma}[Backward parabolic maximal regularity]
\label{lem:forward_parabolic_estimate}
Let $(a,b)\in C^{m,\alpha}(Q;\mathbb E^d)$ and suppose that there exist constants
$0<\lambda_a\leq\Lambda_a<\infty$
such that
$\lambda_a I_d\leq a(t,x)\leq\Lambda_a I_d $
for $(t,x)\in Q.$  For every
$g\in L^2(Q)$ and $\zeta\in H^1(\mathbb T^d)$, the backward problem
\[
-\partial_t\varphi-b\cdot\nabla\varphi-a:\nabla^2\varphi=g,
\qquad
\varphi(1,\cdot)=\zeta,
\]
has a unique solution $\varphi\in \mathcal W:=
H^1(0,1;L^2(\mathbb T^d))
\cap L^2(0,1;H^2(\mathbb T^d))
\hookrightarrow C([0,1];H^1(\mathbb T^d))$, satisfying
\begin{equation}
\label{eq:backward_adjoint_MR_estimate}
\|\varphi\|_{H^1(0,1;L^2)}
+\|\varphi\|_{L^2(0,1;H^2)}
+\|\varphi\|_{C([0,1];H^1)}
\leq
C_{a,b}\bigl(\|g\|_{L^2(Q)}+\|\zeta\|_{H^1}\bigr).
\end{equation}
Here $C_{a,b}$ depends only on $\lambda_a$, $\Lambda_a$ and on upper
bounds for $a$ and $b$ in their stated parabolic H\"older spaces.
\end{lemma}

\begin{proof}
Reverse time and write
$\bar b(s,x):=b(1-s,x)$ and $\bar a(s,x):=a(1-s,x)$.  On
$\mathsf V:=H^1(\mathbb T^d)\hookrightarrow
\mathsf H:=L^2(\mathbb T^d)$, introduce
the forms
\begin{align*}
\mathfrak q_1(s;u,v)
&:=
\int_{\mathbb T^d}\bar a_{ij}(s,x)
\,\partial_ju\,\partial_iv\,\mathrm dx
+\int_{\mathbb T^d}uv\,\mathrm dx,
\\
\mathfrak q_2(s;u,v)
&:=
\int_{\mathbb T^d}\partial_i\bar a_{ij}(s,x)
\,\partial_ju\,v\,\mathrm dx
-\int_{\mathbb T^d}\bar b_j(s,x)\partial_ju\,v\,\mathrm dx
-\int_{\mathbb T^d}uv\,\mathrm dx.
\end{align*}
The spaces $\mathsf V,\mathsf H$ are separable.
The form $\mathfrak q_1$ is symmetric,
uniformly coercive on $\mathsf V$, and
Lipschitz in $s$ on
$\mathsf V\times\mathsf V$, since
$m\geq2$ ensures that the first time
derivative of $a$ is bounded on $Q$.  Because
$\nabla_xa,b\in L^\infty(Q)$, $\mathfrak q_2$ is measurable in $s$
and uniformly bounded on $\mathsf V\times\mathsf H$.  Integration by
parts on the torus identifies the $\mathsf H$-realization associated with
$\mathfrak q_1+\mathfrak q_2$ as
\[
\mathbf A(s)u=-\bar b(s)\cdot\nabla u
-\bar a(s):\nabla^2u.
\]

Since $\mathfrak q_1$ is symmetric and
uniformly coercive, its associated
positive self-adjoint operator
$\mathbf A_1(s)$ has the square-root
property: its square root has domain
$\mathsf V$, and
\[
\|\mathbf A_1(s)^{1/2}u\|_{\mathsf H}^2
=\mathfrak q_1(s;u,u),
\qquad u\in\mathsf V.
\]
In particular, there exist constants
$c,C>0$, independent of $s$, such that
\[
c\|u\|_{\mathsf V}^2
\leq
\|\mathbf A_1(s)^{1/2}u\|_{\mathsf H}^2
\leq
C\|u\|_{\mathsf V}^2.
\]
Moreover, $B(s)=I$ is bounded,
coercive, and self-adjoint.
Thus the hypotheses of
\cite[Theorem~4.2]{ArendtDierLaasriOuhabaz2014}
are satisfied.

Writing $v(s,x)=\varphi(1-s,x)$
and $\bar g(s,x)=g(1-s,x)$,
this theorem gives a unique solution
of
\[
\partial_s v+\mathbf A(s)v=\bar g,
\qquad v(0,\cdot)=\zeta,
\]
such that
\[
v\in H^1(0,1;\mathsf H)
\cap C([0,1];\mathsf V),
\qquad
\mathbf A(\cdot)v(\cdot)\in L^2(Q),
\]
with
\[
\begin{aligned}
&\|v\|_{H^1(0,1;L^2)}
+\|\mathbf A(\cdot)v\|_{L^2(Q)}
+\|v\|_{C([0,1];H^1)}\\
&\qquad\leq C_{a,b}
\bigl(\|\bar g\|_{L^2(Q)}
+\|\zeta\|_{H^1}\bigr).
\end{aligned}
\]

It remains to obtain the
$L^2(0,1;H^2)$ regularity.
For each fixed $s$, the operator
$\mathbf A(s)$ is uniformly elliptic,
and its leading coefficients have
uniformly bounded first spatial
derivatives. The elliptic
$H^2$-regularity theorem
\cite[Section~6.3.1, Theorem~1]{Evans2010},
applied on the torus, gives
\[
\|u\|_{H^2}
\leq C_{a,b}
\bigl(
\|\mathbf A(s)u\|_{L^2}
+\|u\|_{L^2}
\bigr),
\]
whenever $u\in H^1$ and
$\mathbf A(s)u\in L^2$.
The constant is independent of $s$
because the ellipticity and coefficient
bounds are uniform.

Applying this estimate to $u=v(s)$
and integrating in time yields
\[
\|v\|_{L^2(0,1;H^2)}
\leq C_{a,b}
\bigl(
\|\bar g\|_{L^2(Q)}
+\|\zeta\|_{H^1}
\bigr).
\]
Reversing time proves
\eqref{eq:backward_adjoint_MR_estimate}.
\end{proof}

\begin{proof}[Proof of Theorem~\ref{prop:frechet_terminal_map}]
Fix $\psi\in C^{m+2,\alpha}_\diamond(\mathbb T^d)$ and write
$z^\psi=(a^\psi,b^\psi)$. By
Lemma~\ref{lem:HJB_control_differentiability}, the feedback map
$\phi\mapsto z^\phi$ is locally of class $C^{1,1}$ from
$C^{m+2,\alpha}_\diamond(\mathbb T^d)$ into
$C^{m,\alpha}(Q;\mathbb E^d)$.

To establish the smooth dependence of the Fokker--Planck solution
on the coefficients, set
\[
\mathscr X=C^{m,\alpha}(Q),\qquad
\mathscr Z=C^{m,\alpha}(Q;\mathbb E^d),\qquad
\mathscr Y=C^{m-2,\alpha}(Q)\times C^{m,\alpha}(\mathbb T^d).
\]
For $z=(a,b)$ in a sufficiently small open neighborhood of $z^\psi$
in $\mathscr Z$, define
\[
\mathscr F(n,z):=
\left(
\partial_t n+\partial_i(b_i n)-\partial_{ij}(a_{ij}n),
\quad n(0,\cdot)-\mu
\right).
\]
The neighborhood is chosen so that $a$ remains uniformly elliptic.
The parabolic H\"older product estimates show that $\mathscr F$
is a $C^\infty$ map from an open subset of
$\mathscr X\times\mathscr Z$ into $\mathscr Y$.

At $(m^\psi,z^\psi)$, we have $\mathscr F(m^\psi,z^\psi)=0$,
and the partial derivative with respect to $n$ is
\[
D_n\mathscr F(m^\psi,z^\psi)[v]
=
\left(
\partial_t v+\partial_i(b_i^\psi v)-\partial_{ij}(a_{ij}^\psi v),
\quad v(0,\cdot)
\right).
\]
By uniform ellipticity and the periodic Cauchy-problem Schauder theorem
\cite[Chapter~IV, Theorem~5.1]{LadyzhenskayaSolonnikovUraltseva1968},
this operator is an isomorphism from $\mathscr X$ onto $\mathscr Y$.

The Banach-space implicit function theorem gives a locally $C^\infty$
map $z\mapsto m(z)\in\mathscr X$ solving the Fokker--Planck equation
with initial datum $\mu$. For $\phi$ sufficiently close to $\psi$,
uniqueness gives $m^\phi=m(z^\phi)$. Since $z\mapsto m(z)$ is
locally $C^\infty$, its first derivative is locally Lipschitz.
The chain rule and the local $C^{1,1}$ regularity of
$\phi\mapsto z^\phi$ therefore show that $\phi\mapsto m^\phi$ is
locally $C^{1,1}$ with values in $C^{m,\alpha}(Q)$. The terminal
trace is bounded and linear, so $T(\phi)=m^\phi(1,\cdot)$ is locally
$C^{1,1}$ with values in $C^{m,\alpha}(\mathbb T^d)$.

To identify its derivative, differentiate
$\mathscr F(m^\phi,z^\phi)=0$ at $\phi=\psi$ in the direction $h$.
Writing $q^h=Dm^\psi[h]$ and $\dot z^h=(\dot a^h,\dot b^h)$ gives
\[
\partial_t q^h
=
-\partial_i(b_i^\psi q^h+\dot b_i^h m^\psi)
+\partial_{ij}(a_{ij}^\psi q^h+\dot a_{ij}^h m^\psi),
\qquad q^h(0,\cdot)=0.
\]
This is precisely the forward equation in
\eqref{eq:linearized_system}. The Fokker--Planck Schauder
isomorphism established above gives existence and uniqueness of
$q^h$ in $C^{m,\alpha}(Q)$. Together with
Lemma~\ref{lem:HJB_control_differentiability}, this proves existence
and uniqueness for the entire linearized system in the stated spaces.

Taking the terminal trace yields $DT_\psi[h]=q^h(1,\cdot)$.
Integrating the equation for $q^h$ over the torus and using its
zero initial condition gives $\int_{\mathbb T^d}q^h(t,x)\,\mathrm dx=0$
for every $t\in[0,1]$. Thus $DT_\psi[h]\in
C^{m,\alpha}_\diamond(\mathbb T^d)$.
\end{proof}

\subsection{Gramian identity and weak invertibility}
\label{subsec:weak-Gramian}

\begin{lemma}[Gramian identity for the linearized terminal map]
\label{lem:bilinear_form}
Under Assumption~\ref{ass:T_diff}, for every $\psi\in C^{m+2,\alpha}_\diamond(\mathbb T^d)$ and
$h,g\in C^{m+2,\alpha}_\diamond(\mathbb T^d)$,
\begin{equation}
\label{eq:Gramian_identity}
\langle -DT_\psi h,g\rangle
=
\int_Q m^\psi D^2L(a^\psi,b^\psi)
[(\dot a^h,\dot b^h),(\dot a^g,\dot b^g)]
\,\mathrm dx\,\mathrm dt.
\end{equation}
In particular, $-DT_\psi$ is symmetric and nonnegative on
$C^{m+2,\alpha}_\diamond(\mathbb T^d)$.
\end{lemma}

\begin{proof}
Let $(w^h,\dot a^h,\dot b^h,q^h)$ and
$(w^g,\dot a^g,\dot b^g,q^g)$ be the solutions of
\eqref{eq:linearized_system} corresponding to
$h,g\in C^{m+2,\alpha}_\diamond(\mathbb T^d)$.

\noindent
By Theorem~\ref{prop:frechet_terminal_map},
$q^h$ has the regularity required for the following
integration by parts. Set
\[
I(t):=\int_{\mathbb T^d}q^h(t,x)w^g(t,x)\,\mathrm dx.
\]
Using the forward equation for $q^h$ and the backward
equation for $w^g$, we obtain
\begin{align*}
\frac{\mathrm d}{\mathrm dt}I(t)
&=\int_{\mathbb T^d}
\left(
\partial_t q^h\,w^g
+q^h\,\partial_t w^g
\right)\,\mathrm dx\\
&=\int_{\mathbb T^d}
\Big[
(b_i^\psi q^h+\dot b_i^h m^\psi)
\partial_i w^g
+(a_{ij}^\psi q^h+\dot a_{ij}^h m^\psi)
\partial_{ij}w^g\\
&\hspace{3cm}
-q^h\big(
b^\psi\cdot\nabla w^g
+a^\psi:\nabla^2w^g
\big)
\Big]\,\mathrm dx\\
&=\int_{\mathbb T^d}
m^\psi\left(
\dot a^h:\nabla^2w^g
+\dot b^h\cdot\nabla w^g
\right)\,\mathrm dx.
\end{align*}
Here we used periodicity to eliminate boundary terms.

Since $q^h(0,\cdot)=0$ and $w^g(1,\cdot)=g$,
integration over $t\in[0,1]$ gives
\[
\begin{aligned}
\langle DT_\psi[h],g\rangle
&=\int_{\mathbb T^d}q^h(1,x)g(x)\,\mathrm dx\\
&=\int_Q m^\psi
\left(
\dot a^h:\nabla^2w^g
+\dot b^h\cdot\nabla w^g
\right)\,\mathrm dx\,\mathrm dt.
\end{aligned}
\]

Now, taking $(A,\beta)=(\dot a^h,\dot b^h)$
in the linearized optimality condition for
$(\dot a^g,\dot b^g)$ yields
\[
\dot a^h:\nabla^2w^g
+\dot b^h\cdot\nabla w^g
=
-D^2L(a^\psi,b^\psi)
[(\dot a^g,\dot b^g),(\dot a^h,\dot b^h)].
\]
By symmetry of $D^2L$, we conclude that
\[
\langle-DT_\psi[h],g\rangle
=
\int_Q m^\psi D^2L(a^\psi,b^\psi)
[(\dot a^h,\dot b^h),(\dot a^g,\dot b^g)]
\,\mathrm dx\,\mathrm dt,
\]
which proves \eqref{eq:Gramian_identity}.

Symmetry of $-DT_\psi$ follows from symmetry of
$D^2L$. Moreover, taking $g=h$ gives
\[
\langle-DT_\psi[h],h\rangle
=
\int_Q m^\psi D^2L(a^\psi,b^\psi)
[\dot z^h,\dot z^h]\,\mathrm dx\,\mathrm dt
\geq 0,
\]
since $m^\psi\geq0$ and $D^2L$ is positive
definite by Assumption~\ref{ass:T_diff}\textup{(iii)}.
\end{proof}

\begin{remark}
The Gramian identity
\eqref{eq:Gramian_identity} is closely related
to the second-variation representations for
fully nonlinear HJB equations studied in
\cite[Sections~4--5]{
ChassagneuxCrepeyLoeperRuelland2026}.
In that work, positivity of the Hamiltonian
second variation is used to establish local
strong convexity of a finite-dimensional value
map. Here, the corresponding quadratic form
acts on terminal-potential perturbations in an
infinite-dimensional function space.
The parabolic terminal estimate yields
the $H^1_\diamond$ coercivity needed for
invertibility of the terminal map derivative.
\end{remark}

For the weak Gramian analysis, we
consider terminal perturbations
$h\in H^1_\diamond(\mathbb T^d)$.
By Lemma~\ref{lem:forward_parabolic_estimate},
there exists a unique $w^h\in\mathcal W$
satisfying
\begin{equation}
\label{eq:backward_variation_sobolev}
\partial_t w^h
+b^\psi\cdot\nabla w^h
+a^\psi:\nabla^2w^h=0,
\qquad
w^h(1,\cdot)=h.
\end{equation}
For $h\in C^{m+2,\alpha}_\diamond$,
this solution agrees with the $ C^{m+2,\alpha}(Q) $ solution of
\eqref{eq:linearized_system}, by uniqueness.

\begin{lemma}[Parabolic characterization of the terminal $H^1$ norm]
\label{lem:parabolic_terminal_estimate}
Let Assumption~\ref{ass:T_diff} hold. Fix $\psi\in C^{m+2,\alpha}_\diamond(\mathbb T^d)$
and $t_0\in(0,1)$, and set $
Q_{t_0}:=[t_0,1]\times\mathbb T^d $.
Let $h\in H^1_\diamond(\mathbb T^d)$, and let
$w^h\in\mathcal W$ be the unique solution of
\eqref{eq:backward_variation_sobolev}.
Then there exists $C_{\psi,t_0}\geq1$ such that
\begin{equation}
\label{eq:parabolic_H1_equivalence}
C_{\psi,t_0}^{-1}\|h\|_{H^1}
\leq
\|\nabla^2w^h\|_{L^2(Q_{t_0})}
\leq
\|\nabla^2w^h\|_{L^2(Q)}
\leq
C_{\psi,t_0}\|h\|_{H^1}.
\end{equation}
\end{lemma}

\begin{proof}
The homogeneous estimate
\eqref{eq:backward_adjoint_MR_estimate} with
$(a,b)=(a^\psi,b^\psi)$ gives
$
\|\nabla^2w^h\|_{L^2(Q)}
\leq C_\psi\|h\|_{H^1}.
$
The middle inequality in
\eqref{eq:parabolic_H1_equivalence} is immediate.

It remains to prove the first inequality
in \eqref{eq:parabolic_H1_equivalence},
namely that
$\|\nabla^2w^h\|_{L^2(Q_{t_0})}$
controls $\|h\|_{H^1}$.
To this end, define
\[
\overline w^h(t):=
\frac{1}{|\mathbb T^d|}
\int_{\mathbb T^d}w^h(t,x)\,\mathrm dx,
\qquad
\widehat w^h(t,x):=
w^h(t,x)-\overline w^h(t).
\]
The second-order Poincar\'e inequality on
the torus gives
\begin{equation}
\label{eq:centered_H2_estimate}
\|\widehat w^h\|_{L^2(t_0,1;H^2)}
\leq C\|\nabla^2w^h\|_{L^2(Q_{t_0})}.
\end{equation}
The equation for $w^h$ also gives in $L^2$
\[
\begin{aligned}
\partial_t\widehat w^h
={}&-b^\psi\cdot\nabla w^h
-a^\psi:\nabla^2w^h\\
&+\frac{1}{|\mathbb T^d|}
\int_{\mathbb T^d}
\bigl(b^\psi\cdot\nabla w^h
+a^\psi:\nabla^2w^h\bigr)\,\mathrm dx.
\end{aligned}
\]
Each component of $\nabla w^h(t)$ has
zero mean, so Poincar\'e's inequality gives
$
\|\nabla w^h(t)\|_{L^2}
\leq C\|\nabla^2w^h(t)\|_{L^2}.
$
Hence
\begin{equation}
\label{eq:centered_time_estimate}
\|\partial_t\widehat w^h\|_{L^2(t_0,1;L^2)}
\leq C_\psi\|\nabla^2w^h\|_{L^2(Q_{t_0})}.
\end{equation}
The Lions--Magenes trace embedding
\[
L^2(t_0,1;H^2)
\cap H^1(t_0,1;L^2)
\hookrightarrow C([t_0,1];H^1)
\]
\cite[Chapter~I, Theorem~3.1]{LionsMagenes1972},
together with
\eqref{eq:centered_H2_estimate}--
\eqref{eq:centered_time_estimate}, yields
\[
\|\widehat w^h(1)\|_{H^1}
\leq C_{\psi,t_0}
\|\nabla^2w^h\|_{L^2(Q_{t_0})}.
\]
Since $w^h(1)=h$ has zero mean,
$\widehat w^h(1)=h$, proving the first
inequality in
\eqref{eq:parabolic_H1_equivalence}.
\end{proof}

\begin{proposition}[Weak invertibility of the terminal Gramian]
\label{prop:weak_inverse_Gramian}
\label{lem:weak_inverse_Gramian}
Under Assumption~\ref{ass:T_diff}, fix
$\psi\in C^{m+2,\alpha}_\diamond(\mathbb T^d)$.
Then $-DT_\psi$ has a unique continuous extension
\[
(-DT_\psi)^{\mathrm w}:
H^1_\diamond(\mathbb T^d)
\longrightarrow H^{-1}_\diamond(\mathbb T^d),
\]
and this extension is an isomorphism.
Equivalently, for every
$r\in H^{-1}_\diamond(\mathbb T^d)$,
there is a unique
$h\in H^1_\diamond(\mathbb T^d)$
such that $(-DT_\psi)^{\mathrm w}h=r$, and
\begin{equation}
\label{eq:weak_inverse_estimate}
\|h\|_{H^1}\leq C_\psi\|r\|_{H^{-1}}.
\end{equation}
\end{proposition}

\begin{proof}
Fix $t_0:=1/2$ and set
$Q_{t_0}:=[t_0,1]\times\mathbb T^d$.
Compactness of $\J u^\psi(Q)$ and
Assumption~\ref{ass:T_diff}\textup{(ii)}--\textup{(iii)}
give constants
$0<c_{L,\psi}\leq C_{L,\psi}<\infty$
such that
\[
c_{L,\psi}|Z|^2
\leq D^2L(a^\psi,b^\psi)[Z,Z]
\leq C_{L,\psi}|Z|^2,
\qquad Z\in\mathbb E^d,
\]
pointwise on $Q$.
Viewing this Hessian as a self-adjoint operator
on $\mathbb E^d$, set
$
G^\psi:=
\bigl(D^2L(a^\psi,b^\psi)\bigr)^{-1}.
$
Then
\begin{equation}
\label{eq:G_bounds}
\frac{1}{C_{L,\psi}}|Z|^2
\leq Z\boldsymbol\cdot G^\psi[Z]
\leq\frac{1}{c_{L,\psi}}|Z|^2,
\qquad Z\in\mathbb E^d.
\end{equation}

Lemma~\ref{lem:state_system_terminal_map}
and compactness give constants
$0<\underline m_{\psi,t_0}
\leq\overline m_\psi<\infty$ such that
\begin{equation}
\label{eq:m_bounds_weak}
m^\psi\leq\overline m_\psi
\quad\text{on }Q,
\qquad
m^\psi\geq\underline m_{\psi,t_0}
\quad\text{on }Q_{t_0}.
\end{equation}

For $h,g\in H^1_\diamond(\mathbb T^d)$,
let $w^h,w^g\in\mathcal W$ be the unique
solutions of
\eqref{eq:backward_variation_sobolev}
with terminal data $h$ and $g$,
respectively, and define
\begin{equation}
\label{eq:weak_Gramian_form}
\mathfrak B_\psi(h,g)
:=
\int_Qm^\psi\,
\J w^h\boldsymbol\cdot G^\psi[\J w^g]
\,\mathrm dx\,\mathrm dt.
\end{equation}

Since $w^h$ and $w^g$ are periodic,
each component of their spatial gradients
has zero mean. The Poincar\'e inequality
therefore gives
\[
\|\nabla w^h\|_{L^2(Q)}
\leq C\|\nabla^2w^h\|_{L^2(Q)},
\qquad
\|\nabla w^g\|_{L^2(Q)}
\leq C\|\nabla^2w^g\|_{L^2(Q)}.
\]
Therefore, the upper estimate of
Lemma~\ref{lem:parabolic_terminal_estimate},
together with \eqref{eq:G_bounds} and
\eqref{eq:m_bounds_weak}, gives
\[
|\mathfrak B_\psi(h,g)|
\leq C_\psi\|h\|_{H^1}\|g\|_{H^1}.
\]
The lower bounds in \eqref{eq:G_bounds}
and \eqref{eq:m_bounds_weak}, together
with the lower terminal estimate
in Lemma~\ref{lem:parabolic_terminal_estimate},
give
\[
\mathfrak B_\psi(h,h)
\geq
\frac{\underline m_{\psi,t_0}}{C_{L,\psi}}
\|\nabla^2w^h\|_{L^2(Q_{t_0})}^2
\geq c_\psi\|h\|_{H^1}^2.
\]
Thus $\mathfrak B_\psi$ is bounded and
coercive on $H^1_\diamond$.

For $h,g\in C^{m+2,\alpha}_\diamond(\mathbb T^d)$,
the linearized optimality condition gives
$
(\dot a^h,\dot b^h)
=-G^\psi[\J w^h],
$
and Lemma~\ref{lem:bilinear_form} yields
$
\langle -DT_\psi h,g\rangle
=\mathfrak B_\psi(h,g).
$
The bounded form defines
$(-DT_\psi)^{\mathrm w}$ by
\[
\langle (-DT_\psi)^{\mathrm w}h,g
\rangle_{H^{-1}_\diamond,H^1_\diamond}
:=\mathfrak B_\psi(h,g).
\]
For fixed
$h\in C^{m+2,\alpha}_\diamond(\mathbb T^d)$,
density and continuity extend the
displayed identity to all
$g\in H^1_\diamond$.
Hence
$(-DT_\psi)^{\mathrm w}h=-DT_\psi h$
in $H^{-1}_\diamond$.
Since
$C^{m+2,\alpha}_\diamond(\mathbb T^d)$
is dense in $H^1_\diamond$, this proves
that $(-DT_\psi)^{\mathrm w}$ is the
unique continuous extension of $-DT_\psi$
from its stated domain.

Finally, the Lax--Milgram theorem shows
that
\[
(-DT_\psi)^{\mathrm w}:
H^1_\diamond\longrightarrow H^{-1}_\diamond
\]
is an isomorphism and gives
\eqref{eq:weak_inverse_estimate}.
\end{proof}

\begin{remark}[Scope of the weak inverse]
\label{rem:scope-weak-inverse}
Proposition~\ref{prop:weak_inverse_Gramian}
provides an inverse in the Sobolev spaces
$H^{-1}_\diamond\to H^1_\diamond$, but this
is not sufficient for the nonlinear error
analysis. Obtaining quadratic Taylor estimates
in $H^{-1}$ would require additional nonlinear
Sobolev regularity of the terminal map.
This is delicate because the underlying fully
nonlinear HJB equation depends on the spatial
Hessian of its solution, and extending its
solution map to an $H^1$ neighborhood is not
straightforward. We therefore work in H\"older spaces, where the local Lipschitz
continuity of $DT$ gives a quadratic Taylor remainder, and establish
the corresponding strong inverse in
Section~\ref{sec:pseudodifferential-strong-inverse}.
\end{remark}

\section{Pseudodifferential analysis and strong invertibility}
\label{sec:pseudodifferential-strong-inverse}

The weak invertibility established in
Proposition~\ref{prop:weak_inverse_Gramian}
does not directly yield an inverse in strong
H\"older spaces.
To obtain such an inverse, we identify the
terminal Gramian as an elliptic classical
pseudodifferential operator of order two.
We first recall the necessary notation.

\subsection{Pseudodifferential notation}
\label{subsec:pseudodifferential-notation}

We recall the notation for classical pseudodifferential
operators on $\mathbb T^d$.

Let $\Psi^{-\infty}(\mathbb T^d)$ denote the class
of smoothing operators, namely operators with
smooth Schwartz kernels on
$\mathbb T^d\times\mathbb T^d$. For $r\in\mathbb R$, a smooth function
$p:\mathbb T^d\times\mathbb R^d\to\mathbb C$
belongs to $S^r$ if, for all multi-indices
$\gamma,\delta$, there exists $C_{\gamma,\delta}<\infty$
such that
\[
|\partial_x^\gamma\partial_\xi^\delta
p(x,\xi)|
\leq C_{\gamma,\delta}
(1+|\xi|)^{r-|\delta|}.
\]

\noindent
The periodic quantization of $p$ is
\[
\operatorname{Op}(p)u(x)
:=
\sum_{n\in\mathbb Z^d}
e^{2\pi i n\cdot x}
p(x,2\pi n)\widehat u(n),
\]
where $\widehat u(n)$ denotes the Fourier
coefficient of $u$.

\noindent
The class $\Psi^r(\mathbb T^d)$ consists of
operators of the form
\[
A=\operatorname{Op}(p)+R,
\qquad
p\in S^r,\quad
R\in\Psi^{-\infty}(\mathbb T^d).
\]
With this definition,
\[
\Psi^{-\infty}(\mathbb T^d)
=
\bigcap_{r\in\mathbb R}\Psi^r(\mathbb T^d).
\]

\noindent
A symbol $p\in S^r$ is called classical,
written $p\in S^r_{\mathrm{cl}}$, if it admits
an asymptotic expansion
\[
p\sim\sum_{j=0}^{\infty}p_{r-j},
\]
where $p_{r-j}(x,\xi)$ is positively homogeneous
of degree $r-j$ in $\xi$ for $|\xi|\geq1$, and the asymptotic representation means that
\[
p-\sum_{j=0}^{N-1}p_{r-j}
\in S^{r-N},
\qquad N\geq1.
\]
Here the homogeneous terms are understood
to be smoothly extended near $\xi=0$.
The corresponding operator class is denoted
$\Psi^r_{\mathrm{cl}}(\mathbb T^d)$.
For $A\in\Psi^r_{\mathrm{cl}}$, the leading
homogeneous term $p_r$ is its principal symbol,
denoted by $\sigma_r(A)$.

An operator $A\in\Psi^r_{\mathrm{cl}}$ is called an elliptic classical pseudodifferential
operator of order $r$
if there exist $c,R>0$ such that
\[
|\sigma_r(A)(x,\xi)|
\geq c|\xi|^r,
\qquad |\xi|\geq R.
\]

\noindent
We use the standard composition, adjoint,
asymptotic summation, and elliptic parametrix
results of the classical pseudodifferential
calculus; see
\cite{Hormander1994III,Taylor1991Pseudo}.

For the decomposition
$\mathbb E^d=\mathbb S^d\times\mathbb R^d$,
we denote by
\[
\iota_a:\mathbb S^d\longrightarrow\mathbb E^d,
\qquad
\iota_a(N):=(N,0),
\]
the canonical inclusion, and by
\[
\pi_a:\mathbb E^d\longrightarrow\mathbb S^d,
\qquad
\pi_a(A,b):=A,
\]
the projection onto the diffusion component.

\subsection{Regularity on a terminal collar}
\label{subsec:terminal-collar}

\begin{lemma}[Smoothness on a terminal collar]
\label{lem:smooth-terminal-coefficients}
Let Assumption~\ref{ass:T_diff} hold, and
suppose in addition that
$
L\in C^\infty_{\mathrm{loc}}(\operatorname{int}U)$,
$
\psi\in C^\infty_\diamond(\mathbb T^d).
$
Set
$
c^\psi
:=
m^\psi
\bigl(D^2L(a^\psi,b^\psi)\bigr)^{-1}.
$
Then there exists $\delta_\psi\in(0,1/2)$
such that
\[
a^\psi,b^\psi,c^\psi
\in C^\infty(
[1-2\delta_\psi,1]\times\mathbb T^d).
\]
\end{lemma}

\begin{proof}
The optimizer satisfies
\[
 DL(\mathfrak z(M,p))+(M,p)=0.
\]
Since it lies in the interior of \(U\) and
\(D^2L(\mathfrak z(M,p))\) is invertible, the implicit-function theorem
at every finite order gives
$
 \mathfrak z\in C^\infty_{\mathrm{loc}}
 (\mathbb E^d;\operatorname{int}U).
$
The envelope identity \(DH=\mathfrak z\) therefore implies
\(H\in C^\infty_{\mathrm{loc}}(\mathbb E^d)\).

Let \(v(s,x):=u^\psi(1-s,x)\).  It solves
\[
 \partial_s v-H(\nabla^2v,\nabla v)=0,
 \qquad v(0,\cdot)=\psi.
\]
The compactness of \(\J v(Q)\) bounds every derivative of \(H\) near
this range.  Each spatial derivative of \(v\) solves a linear uniformly
parabolic equation with principal and first-order coefficients
\[
 D_MH(\J v)=\mathfrak a(\J v),
 \qquad D_pH(\J v)=\mathfrak b(\J v),
\]
and lower-order terms involving only previously controlled derivatives.
The periodic Cauchy-problem Schauder estimates
\cite[Chapter~IV, Theorem~5.1 and Section~14]{LadyzhenskayaSolonnikovUraltseva1968}
and induction control all spatial derivatives; uniform ellipticity follows
from Assumption~\ref{ass:T_diff}\textup{(iv)}.  The equation then yields
all time derivatives, so
$
 v\in C^\infty([0,1]\times\mathbb T^d).
$
Consequently, $(a^\psi,b^\psi)$ and
$\bigl(D^2L(a^\psi,b^\psi)\bigr)^{-1}$
are smooth on $Q$.

Expanding the Fokker--Planck equation gives
\[
 \partial_t m^\psi
 =a^\psi_{ij}\partial_{ij}m^\psi
  +\bigl(2\partial_i a^\psi_{ij}-b^\psi_j\bigr)\partial_jm^\psi
  +\bigl(\partial_{ij}a^\psi_{ij}
  -\partial_i b^\psi_i\bigr)m^\psi.
\]
This equation is uniformly parabolic with smooth coefficients.  Interior
Schauder estimates on $[\varepsilon/2,1]\times\mathbb T^d$ and
bootstrapping (equivalently, the positive-time linear theory in
\cite[Chapter~IV, Sections~10 and~14]{LadyzhenskayaSolonnikovUraltseva1968})
give
\[
 m^\psi\in
 C^\infty([\varepsilon,1]\times\mathbb T^d)
 \qquad\text{for every }\varepsilon>0,
\]
without additional regularity of \(\mu\) beyond
Assumption~\ref{ass:T_diff}\textup{(v)}.  Thus
$c^\psi
=m^\psi\bigl(D^2L(a^\psi,b^\psi)\bigr)^{-1}$
is smooth on every terminal collar,
and any $\delta_\psi\in(0,1/2)$
suffices.
\end{proof}

\subsection{Pseudodifferential structure of the terminal Gramian}
\label{subsec:smooth-terminal-Gramian}

We establish a pseudodifferential
representation for a parabolic Gramian
whose coefficients are smooth near the
initial time $s=0$, but have only finite
regularity elsewhere.
In our application, these coefficients
arise by reversing time in the original
HJB--Fokker--Planck system, so the required
smoothness near $s=0$ follows from
Lemma~\ref{lem:smooth-terminal-coefficients}.
The principal symbol depends only on the
coefficients at $s=0$, while the contribution
from times bounded away from zero is smoothing.

\begin{lemma}[Pseudodifferential structure and ellipticity of the Gramian]
\label{lem:smooth-terminal-Gramian}
Let
$(\bar a,\bar b)
\in C^{m,\alpha}(Q;\mathbb E^d) $,
$\bar c\in
C^{m,\alpha}\bigl(
Q;\mathcal L(\mathbb E^d;\mathbb E^d)
\bigr)
$,
and assume that $\bar c(s,x)$ is
self-adjoint. Suppose that there exist
constants $0<\lambda\leq\Lambda<\infty$
such that
\[
\lambda|\xi|^2
\leq
\bar a(s,x):(\xi\otimes\xi)
\leq
\Lambda|\xi|^2
\]
for every $(s,x)\in Q$ and
$\xi\in\mathbb R^d$.

Assume also that there exist \(\delta\in(0,1/2)\) and \(\kappa>0\)
such that
$
 \bar a,\bar b,\bar c\in
 C^\infty([0,2\delta]\times\mathbb T^d)
$
and
$Z\boldsymbol\cdot\bar c(0,x)[Z]\geq\kappa|Z|^2$ for $x\in\mathbb T^d$ and $\ Z\in\mathbb E^d$. For every $h\in H^1(\mathbb T^d)$,
let $P_\cdot h\in\mathcal W$ be the unique
solution of
\[
\partial_sP_sh
-\bar b(s,x)\cdot\nabla P_sh
-\bar a(s,x):\nabla^2P_sh=0,
\qquad
P_0h=h,
\]
and define \(A:H^1\to H^{-1}\) by
\begin{equation}
 \langle Ah,g\rangle
 :=\int_0^1\!\int_{\mathbb T^d}
 \J P_sh\boldsymbol\cdot\bar c(s,x)[\J P_sg]
 \,\mathrm dx\,\mathrm ds.
 \label{eq:smooth-collar-Gramian-form}
\end{equation}
Then \(A\) is the restriction to \(H^1\) of a classical
pseudodifferential operator in \(\Psi^2_{\mathrm{cl}}(\mathbb T^d)\).
Its principal symbol is
\begin{equation}
 \sigma_2(A)(x,\xi)
 =\frac{
  \iota_a(\xi\otimes\xi)\boldsymbol\cdot
  \bar c(0,x)[\iota_a(\xi\otimes\xi)]
 }{
  2\,\bar a(0,x):(\xi\otimes\xi)
 },
 \qquad \xi\neq0.
 \label{eq:smooth-collar-principal-symbol}
\end{equation}
In particular, there exist \(c,C>0\) such that
\[
 c|\xi|^2\leq\sigma_2(A)(x,\xi)\leq C|\xi|^2
 \qquad(\xi\neq0),
\]
so \(A\) is an elliptic classical
pseudodifferential operator of order two.
\end{lemma}

\begin{proof}
Choose \(\chi\in C_c^\infty([0,2\delta))\), equal to one on
\([0,\delta]\), and split
$
A=A_{\mathrm{near}}+A_{\mathrm{far}}
$
by inserting \(\chi(s)\) and \(1-\chi(s)\), respectively, in
\eqref{eq:smooth-collar-Gramian-form}.

\medskip
\noindent\emph{The far-time contribution.}

Let $G(s,t)$, $0\leq t\leq s\leq1$,
denote the evolution operator associated
with
\[
B_s:=\bar b(s,x)\cdot\nabla
+\bar a(s,x):\nabla^2.
\]
Thus $G(s,t)v$ solves
$\partial_su=B_su$ with $u(t)=v$,
and $P_s=G(s,0)$.

By Lemma~\ref{lem:forward_parabolic_estimate},
applied after time reversal,
$P_\cdot h\in\mathcal W$ for every
$h\in H^1(\mathbb T^d)$.
Since $\bar c$ is bounded, the
Gramian form
\eqref{eq:smooth-collar-Gramian-form}
is therefore well defined.

Set $s_*:=\delta/2$.
Applying Lemma~\ref{lem:forward_parabolic_estimate},
after time reversal, to the evolution on
$[s_*,1]$ gives
\begin{equation}
\label{eq:rough-evolution-H1-H2}
\sup_{s_*\leq s\leq1}
\|G(s,s_*)v\|_{H^1}^2
+
\int_{s_*}^1
\|G(s,s_*)v\|_{H^2}^2\,\mathrm ds
\leq C\|v\|_{H^1}^2.
\end{equation}
Since $s\geq\delta>s_*$ on
$\operatorname{supp}(1-\chi)$, the
evolution property gives
$
P_s=G(s,s_*)P_{s_*}.
$
Define
\[
\begin{aligned}
\mathfrak m(v,w):={}&
\int_{s_*}^1(1-\chi(s))
\int_{\mathbb T^d}
\J G(s,s_*)v
\boldsymbol\cdot\bar c(s)
[\J G(s,s_*)w]
\,\mathrm dx\,\mathrm ds.
\end{aligned}
\]
By Cauchy--Schwarz and
\eqref{eq:rough-evolution-H1-H2},
$
|\mathfrak m(v,w)|
\leq C\|\bar c\|_{L^\infty}
\|v\|_{H^1}\|w\|_{H^1}.
$
Hence $\mathfrak m$ defines
$M\in\mathcal L(H^1,H^{-1})$, and
$
A_{\mathrm{far}}=P_{s_*}^*MP_{s_*}.
$

Since the coefficients are smooth
on $[0,s_*]$ and uniformly elliptic,
positive-time parabolic regularization
gives
$
P_{s_*},P_{s_*}^*:
\mathcal D'(\mathbb T^d)
\longrightarrow C^\infty(\mathbb T^d).
$
Consequently, $P_{s_*}^*MP_{s_*}$ has a smooth
Schwartz kernel, and therefore
\begin{equation}
\label{eq:far-smoothing}
A_{\mathrm{far}}
\in\Psi^{-\infty}(\mathbb T^d).
\end{equation}

\medskip
\noindent\emph{The near-time contribution.}
Set $s_1:=2\delta$.
Recall that $\partial_sP_s=B_sP_s$ and $P_0=I$,
where $B_s$ was defined above.

Define the fourth-order differential operator
$
Q_s v
:=
\chi(s)\mathcal J^*
\bigl(\bar c(s,x)[\mathcal Jv]\bigr),
0\leq s\leq s_1,
$
where $\mathcal J^*$ denotes the formal
$L^2$-adjoint of $\mathcal J$.
By integration by parts on the torus,
\[
\langle Q_s v,u\rangle
=
\chi(s)\int_{\mathbb T^d}
\mathcal Ju\boldsymbol\cdot
\bar c(s,x)[\mathcal Jv]
\,\mathrm dx.
\]
Consequently,
\[
A_{\mathrm{near}}
=
\int_0^{s_1}P_s^*Q_sP_s\,\mathrm ds.
\]

We construct a smooth parameter-dependent family
$
C_s\in\Psi^2_{\mathrm{cl}}(\mathbb T^d),
0\leq s\leq s_1,
$
such that
\begin{equation}
\label{eq:symbolic-Lyapunov-remainder}
R_s
:=
Q_s-\bigl(
\partial_sC_s+B_s^*C_s+C_sB_s
\bigr)
\in\Psi^{-\infty}(\mathbb T^d)
\end{equation}
smoothly and uniformly in $s$.

Write
$
l_s(x,\xi)
:=
\bar a(s,x):(\xi\otimes\xi).
$
Then
$
\lambda|\xi|^2\leq l_s(x,\xi)\leq\Lambda|\xi|^2.
$
For a smooth family
$I_s\in\Psi^r_{\mathrm{cl}}$, set
$
\mathcal L_s(I_s)
:=
\partial_sI_s+B_s^*I_s+I_sB_s.
$
The classical adjoint and composition formulas give
\begin{equation}
\label{eq:Lyapunov-leading-symbol}
\mathcal L_s(I_s)\in\Psi^{r+2}_{\mathrm{cl}},
\qquad
\sigma_{r+2}\bigl(\mathcal L_s(I_s)\bigr)
=-2l_s\,\sigma_r(I_s),
\end{equation}
because both $B_s$ and $B_s^*$ have principal
symbol $-l_s$.

These formulas and the asymptotic summation below
are the parameter-dependent versions of
\cite[Proposition~18.1.3 and
Theorems~18.1.7--18.1.8]{Hormander1994III};
coordinate invariance and globalization follow from
\cite[Theorem~18.1.17]{Hormander1994III}.
Quantize the recursively constructed homogeneous
symbols in torus charts and patch them by a partition
of unity. Since $s\in[0,s_1]$, choosing the
asymptotic-summation cutoffs independent of $s$
and the excision radii diagonally over the joint
seminorms in $(s,x,\xi)$ preserves the symbol
estimates after every $s$-derivative.

The principal symbol of $Q_s$ is
\[
q_4(s,x,\xi)
=
\chi(s)\,
\iota_a(\xi\otimes\xi)\boldsymbol\cdot
\bar c(s,x)[\iota_a(\xi\otimes\xi)].
\]
The drift--drift and cross blocks of $\bar c$
have differential orders two and at most three,
respectively; only the diffusion--diffusion block
has order four.

Choose $C_s^{(0)}\in\Psi^2_{\mathrm{cl}}$
with principal symbol
\[
c_2(s,x,\xi)
=-\frac{q_4(s,x,\xi)}{2l_s(x,\xi)}.
\]
Then
$Q_s-\mathcal L_s(C_s^{(0)})
\in\Psi^3_{\mathrm{cl}}$.
Inductively, for $N\geq1$, if terms through
degree $3-N$ give a residual of order at most
$4-N$ with homogeneous principal symbol
$r_{4-N}$, add a term of order $2-N$
with homogeneous principal symbol
\[
c_{2-N}:=-\frac{r_{4-N}}{2l_s}.
\]
Define this quotient for $|\xi|\geq1$
and extend it smoothly across the zero section
before quantization, without changing its
homogeneous principal symbol.
Equation~\eqref{eq:Lyapunov-leading-symbol}
lowers the residual by one order.
Since $\partial_sC_s$ lies two orders below
the leading term of $\mathcal L_s(C_s)$,
it involves only previously constructed terms
at each stage.
Smooth parameter-dependent asymptotic summation
now gives $C_s$ satisfying
\eqref{eq:symbolic-Lyapunov-remainder}.

For smooth initial data
$h,g\in C^\infty(\mathbb T^d)$,
the regularity of the coefficients on
$[0,s_1]$ justifies differentiating
$P_s^*C_sP_s$. Using
$\partial_sP_s=B_sP_s$ and
\eqref{eq:symbolic-Lyapunov-remainder},
we obtain
\[
\begin{aligned}
\frac{\mathrm d}{\mathrm ds}
\bigl(P_s^*C_sP_s\bigr)
&=
P_s^*
\bigl(
\partial_sC_s+B_s^*C_s+C_sB_s
\bigr)P_s\\
&=P_s^*(Q_s-R_s)P_s.
\end{aligned}
\]
Integrating over $[0,s_1]$ and using
$P_0=I$ yields
\begin{equation}
\label{eq:near-Gramian-mod-smoothing}
A_{\mathrm{near}}
=
P_{s_1}^*C_{s_1}P_{s_1}-C_0
+\int_0^{s_1}P_s^*R_sP_s\,\mathrm ds.
\end{equation}

The term $P_{s_1}^*C_{s_1}P_{s_1}$
is smoothing because $s_1>0$ and
$P_{s_1}$ has a smooth kernel.
Differentiated energy estimates bound $P_s$
uniformly on every nonnegative integer Sobolev
space for $0\leq s\leq s_1$;
the time-reversed adjoint equation gives the same
bounds for $P_s^*$, and duality and interpolation
extend them to every real order.

For arbitrary $M,N\geq0$, the uniformly
smoothing $R_s$ maps $H^{-N}$ to $H^M$
uniformly. The same therefore holds for
$P_s^*R_sP_s$ and its integral, whose
Schwartz kernel is smooth.

It follows from
\eqref{eq:near-Gramian-mod-smoothing} that
$
A_{\mathrm{near}}=-C_0
\quad\bmod\Psi^{-\infty}(\mathbb T^d).
$
Since $\chi(0)=1$, its principal symbol is
\[
\sigma_2(A_{\mathrm{near}})(x,\xi)
=
\frac{
\iota_a(\xi\otimes\xi)\boldsymbol\cdot
\bar c(0,x)[\iota_a(\xi\otimes\xi)]
}{
2\,\bar a(0,x):(\xi\otimes\xi)
}.
\]
In particular,
\[
\frac{\kappa}{2\Lambda}|\xi|^2
\leq
\sigma_2(A_{\mathrm{near}})(x,\xi)
\leq
\frac{
\sup_{x\in\mathbb T^d}
\|\bar c(0,x)\|_{\mathcal L(\mathbb E^d;\mathbb E^d)}
}{2\lambda}
|\xi|^2.
\]
Together with \eqref{eq:far-smoothing},
this identifies the Gramian on smooth
functions with an elliptic classical
pseudodifferential operator of order two.
Since both this operator and the
form-defined operator $A$ are continuous
from $H^1$ to $H^{-1}$, density extends
the identification to $H^1$.
\end{proof}

\subsection{Proof of the strong inverse result}
\label{subsec:proof-strong-inverse}

\begin{proof}[Proof of Theorem~\ref{prop:strong-inverse-smooth-terminal}]
After time reversal, Lemma~\ref{lem:smooth-terminal-coefficients}
verifies the smooth-collar hypotheses of
Lemma~\ref{lem:smooth-terminal-Gramian} for
\[
 \begin{aligned}
 \bar a(s,x)&=a^\psi(1-s,x),\\
 \bar b(s,x)&=b^\psi(1-s,x),\\
 \bar c(s,x)&=m^\psi(1-s,x)
 \bigl(D^2L(a^\psi(1-s,x),b^\psi(1-s,x))\bigr)^{-1}.
 \end{aligned}
\]
By Lemma~\ref{lem:state_system_terminal_map}
and Assumption~\ref{ass:T_diff},
the coefficients $\bar a,\bar b,\bar c$
belong to $C^{m,\alpha}(Q)$ and satisfy
the uniform ellipticity condition of
Lemma~\ref{lem:smooth-terminal-Gramian}.
Their smoothness near $s=0$ follows
from Lemma~\ref{lem:smooth-terminal-coefficients}.

Moreover, by
\eqref{eq:strict_positivity_density},
$m^\psi(1,x)>0$ for every
$x\in\mathbb T^d$. Continuity and
compactness therefore imply
$
\min_{x\in\mathbb T^d}m^\psi(1,x)>0.
$
By \eqref{eq:G_bounds}, for every $Z\in\mathbb E^d$,
\[
Z\boldsymbol\cdot\bar c(0,x)[Z]
\geq \frac{\min_{y\in\mathbb T^d}m^\psi(1,y)}{C_{L,\psi}}\,|Z|^2.
\]
Thus $\bar c(0,x)$ is uniformly positive definite.
Hence Lemma~\ref{lem:smooth-terminal-Gramian}
implies that $-DT_\psi$ is an elliptic
classical pseudodifferential operator
of order two.

Set
\[
G^\psi
:=
\bigl(D^2L(a^\psi,b^\psi)\bigr)^{-1},
\qquad
G^\psi_{aa}:=\pi_aG^\psi\iota_a.
\]
By \eqref{eq:smooth-collar-principal-symbol},
\[
\sigma_2(-DT_\psi)(x,\xi)
=
\frac{
m^\psi(1,x)(\xi\otimes\xi):
G^\psi_{aa}(1,x)[\xi\otimes\xi]
}{
2a^\psi(1,x):(\xi\otimes\xi)
},
\qquad \xi\neq0.
\]

The classical elliptic parametrix theorem on \(\mathbb T^d\) (see \cite[Section~0.4]{Taylor1991Pseudo}) gives
\(E_\psi\in\Psi^{-2}_{\mathrm{cl}}\) and
\(R_\psi\in\Psi^{-\infty}\) such that
\begin{equation}
 E_\psi (-DT_\psi)=I-R_\psi.
 \label{eq:classical-left-parametrix}
\end{equation}

\noindent
The H\"older--Zygmund mapping theorem gives
\[
E_\psi:
C^{m,\alpha}(\mathbb T^d)
\longrightarrow
C^{m+2,\alpha}(\mathbb T^d),
\qquad
R_\psi:
H^1(\mathbb T^d)
\longrightarrow C^\infty(\mathbb T^d).
\]
See \cite[Appendix~A.1, especially
(A.1.6)--(A.1.7)]{Taylor1991Pseudo}.

Let
$r\in C^{m,\alpha}_\diamond(\mathbb T^d)$.
By Proposition~\ref{prop:weak_inverse_Gramian},
there is a unique
$h\in H^1_\diamond(\mathbb T^d)$
such that
$(-DT_\psi)^{\mathrm w}h=r$, with
\[
\|h\|_{H^1}
\leq C_\psi\|r\|_{H^{-1}}.
\]

After the time reversal $s=1-t$, the
Gramian form in
\eqref{eq:smooth-collar-Gramian-form}
coincides with $\mathfrak B_\psi$ in
\eqref{eq:weak_Gramian_form}.
Consequently, the pseudodifferential
realization of
Lemma~\ref{lem:smooth-terminal-Gramian}
agrees on $H^1_\diamond$ with the weak
extension $(-DT_\psi)^{\mathrm w}$ of
Proposition~\ref{prop:weak_inverse_Gramian}. Applying
\eqref{eq:classical-left-parametrix}
to the weak equation in
$\mathcal D'(\mathbb T^d)$ gives
$
h=E_\psi r+R_\psi h
\in C^{m+2,\alpha}(\mathbb T^d).
$
Moreover,
\[
\begin{aligned}
\|h\|_{C^{m+2,\alpha}}
&\leq C_\psi
\bigl(
\|r\|_{C^{m,\alpha}}
+\|h\|_{H^1}
\bigr)\\
&\leq C'_\psi\|r\|_{C^{m,\alpha}},
\end{aligned}
\]
where we used the continuous embedding
$C^{m,\alpha}\hookrightarrow H^{-1}$. Thus $-DT_\psi$ is surjective from
$C^{m+2,\alpha}_\diamond$ onto
$C^{m,\alpha}_\diamond$ and has a bounded
inverse. Injectivity follows from
Proposition~\ref{prop:weak_inverse_Gramian}.
Hence
\[
DT_\psi:
C^{m+2,\alpha}_\diamond
\xrightarrow{\ \sim\ }
C^{m,\alpha}_\diamond.
\]

Finally, if $r\in C^\infty_\diamond$,
the identity $h=E_\psi r+R_\psi h$
shows that $h\in C^\infty_\diamond$,
since pseudodifferential operators preserve
smoothness and $R_\psi$ is smoothing.
Thus the inverse also preserves smoothness.
\end{proof}

\section{Proofs of the linearized convergence results}
\label{sec:proofs-linearized}

\subsection{Proof of Corollary~\ref{prop:smooth-induction-geometric-directions}}

\begin{proof}
We proceed by induction on $k$.
By assumption,
$\psi_0\in C^\infty_\diamond$.
Suppose that
$\psi_k\in C^\infty_\diamond$.
Since
$L\in C^\infty_{\mathrm{loc}}(\operatorname{int}U)$,
Theorem~\ref{prop:strong-inverse-smooth-terminal}
implies that
\[
-DT_{\psi_k}:
C^{m+2,\alpha}_\diamond
\xrightarrow{\ \sim\ }
C^{m,\alpha}_\diamond
\]
is an isomorphism whose inverse
preserves smoothness. Since
$
d_k=\tau(1+\tau)^{-k-1}\rho_0
\in C^\infty_\diamond,
$
there exists a unique
$
h_k=(-DT_{\psi_k})^{-1}d_k
\in C^\infty_\diamond.
$
Consequently,
$
\psi_{k+1}=\psi_k+h_k
\in C^\infty_\diamond.
$
This closes the induction and proves that the
recursion is well defined for every $k\geq0$,
with smooth iterates and increments.
\end{proof}

\subsection{Proof of Theorem~\ref{thm:open_loop_error}}

Throughout the remaining proofs, write
\[
\mathsf X:=C^{m+2,\alpha}_\diamond(\mathbb T^d),
\qquad
\mathsf Y:=C^{m,\alpha}_\diamond(\mathbb T^d).
\]
The set $\calK_R$ denotes the closed
$\mathsf X$-ball of radius $R$ centered at $\psi_0$.

\begin{lemma}[Local uniform invertibility of the linearized terminal map]
\label{prop:local-uniform-inverse}
Under Assumption~\ref{ass:open-loop-error},
there exist finite, strictly positive $R_0$, $M$ and $L_A$
such that, for every $\psi\in\calK_{R_0}$,
$
-DT_\psi:
\mathsf X\xrightarrow{\ \sim\ }\mathsf Y
$
is an isomorphism satisfying
\[
\|(-DT_\psi)^{-1}\|_{\mathcal L(\mathsf Y,\mathsf X)}
\leq M.
\]
Moreover,
\begin{equation}
\label{eq:A-local-operator-Lipschitz}
\|(-DT_\psi)-(-DT_\phi)\|_{\mathcal L(\mathsf X,\mathsf Y)}
\leq L_A\|\psi-\phi\|_{\mathsf X},
\qquad
\psi,\phi\in\calK_{R_0}.
\end{equation}
\end{lemma}

\begin{proof}
By Theorem~\ref{prop:frechet_terminal_map},
the terminal map $T$ is locally $C^{1,1}$
from $\mathsf X$ into $C^{m,\alpha}$.
\noindent
Since $\psi_0\in C^\infty_\diamond$ and
$L\in C^\infty_{\mathrm{loc}}(\operatorname{int}U)$,
Theorem~\ref{prop:strong-inverse-smooth-terminal}
implies that
$
(-DT_{\psi_0}):\mathsf X\xrightarrow{\ \sim\ }\mathsf Y.
$
Set
\[
M_0 := \left\|(-DT_{\psi_0})^{-1}\right\|_{\mathcal{L}(\mathsf{Y},\mathsf{X})},
\qquad
M := 2M_0.
\]

\noindent
Since $T$ is locally $C^{1,1}$ at $\psi_0$ and $DT_\psi$
takes values in $\mathsf Y$, there exist $R_1>0$ and
$L_A<\infty$ such that
$
\|(-DT_\psi)-(-DT_\phi)\|_{\mathcal L(\mathsf X,\mathsf Y)}
\leq L_A\|\psi-\phi\|_{\mathsf X},
\psi,\phi\in\calK_{R_1}.
$
\noindent
Choose $R_0\in(0,R_1]$ sufficiently small that
$
M_0L_AR_0\leq\frac12.
$
Then, for every $\psi\in\calK_{R_0}$,
$
\|(-DT_{\psi_0})^{-1}((-DT_\psi)-(-DT_{\psi_0}))\|
_{\mathcal L(\mathsf X)}
\leq M_0L_AR_0\leq\frac12.
$
The factorization
\[
(-DT_\psi)
=
(-DT_{\psi_0})
\left[
I+(-DT_{\psi_0})^{-1}((-DT_\psi)-(-DT_{\psi_0}))
\right]
\]
and the Neumann-series argument show that
$(-DT_\psi)$ is invertible, with
$
\|(-DT_\psi)^{-1}\|_{\mathcal L(\mathsf Y,\mathsf X)}
\leq 2M_0=M.
$
Finally, the previously established Lipschitz
estimate on $\calK_{R_1}$ remains valid on
$\calK_{R_0}$, proving
\eqref{eq:A-local-operator-Lipschitz}.
\end{proof}

\begin{lemma}[Confinement of the open-loop iterates]
\label{lem:iterate-confinement}
Under Assumption~\ref{ass:open-loop-error},
let $R_0,M$ be supplied by
Lemma~\ref{prop:local-uniform-inverse}.
Fix $R\in(0,R_0]$ and $\tau>0$.
If
$
M\|\rho_0\|_{\mathsf Y}\leq R,
$
then \eqref{eq:recursion} is uniquely defined for every $k\geq0$,
with $\psi_k\in\calK_R$ and
\begin{equation}
\|\psi_k-\psi_0\|_{\mathsf X}
\leq M\bigl(1-(1+\tau)^{-k}\bigr)\|\rho_0\|_{\mathsf Y}<R.
\label{eq:iterate-local-confinement}
\end{equation}
Moreover,
\begin{equation}
\|h_k\|_{\mathsf X}\leq M\tau(1+\tau)^{-k-1}\|\rho_0\|_{\mathsf Y},
\label{eq:local-geometric-h-bound}
\end{equation}
and
\begin{equation}
\sup_{k\geq0}\|(-DT_{\psi_k})^{-1}\|_{\mathcal L(\mathsf Y,\mathsf X)}\leq M.
\label{eq:uniform-in-k-strong-inverse}
\end{equation}
\end{lemma}

\begin{proof}
We construct the iterates by induction, using the uniform inverse
bound on the whole ball $\calK_R$. Suppose $\psi_0,\ldots,\psi_k$
have been constructed in that ball. For
$d_k=\tau(1+\tau)^{-k-1}\rho_0\in\mathsf Y$, by Lemma \ref{prop:local-uniform-inverse} the equation
$(-DT_{\psi_k})h_k=d_k$ has a unique solution in $\mathsf X$, and
\eqref{eq:local-geometric-h-bound} follows. Defining
$\psi_{k+1}:=\psi_k+h_k$, we obtain
\[
\begin{aligned}
\|\psi_{k+1}-\psi_0\|_{\mathsf X}
&\leq\sum_{j=0}^{k}\|h_j\|_{\mathsf X}\\
&\leq M\|\rho_0\|_{\mathsf Y}\sum_{j=0}^{k}\tau(1+\tau)^{-j-1}\\
&=M\bigl(1-(1+\tau)^{-k-1}\bigr)\|\rho_0\|_{\mathsf Y}<R.
\end{aligned}
\]
Thus $\psi_{k+1}\in\calK_R$, which closes the induction and proves
\eqref{eq:iterate-local-confinement}. The uniform inverse bound of Lemma \ref{prop:local-uniform-inverse}
then gives \eqref{eq:uniform-in-k-strong-inverse}.
\end{proof}

\begin{proof}[Proof of Theorem~\ref{thm:open_loop_error}]
Let $R_0,M,L_A$ be the constants supplied by
Lemma~\ref{prop:local-uniform-inverse}. Fix $R\in(0,R_0]$
and $\tau>0$ such that $M\|\rho_0\|_{\mathsf Y}\leq R$.
By Lemma~\ref{lem:iterate-confinement}, the recursion is
uniquely defined for every $k\geq0$, its iterates remain
in $\calK_R$, and
$\|(-DT_{\psi_k})^{-1}\|_{\mathcal L(\mathsf Y,\mathsf X)}
\leq M$ uniformly in $k$. 

\noindent
Define the one-step Taylor remainder by
$e_k:=T(\psi_{k+1})-T(\psi_k)-DT_{\psi_k}[h_k]$.
Since $\calK_R$ is convex and contains both endpoints
of the update, the segment $\psi_k+sh_k$, $s\in[0,1]$,
lies in $\calK_R$. The integral Taylor formula therefore gives
\[
e_k=\int_0^1
\bigl(DT_{\psi_k+sh_k}-DT_{\psi_k}\bigr)[h_k]
\,\mathrm ds.
\]
Using \eqref{eq:A-local-operator-Lipschitz} and
\eqref{eq:local-geometric-h-bound}, we obtain
\begin{equation}
\label{eq:error-geometric-one-step-bound}
\|e_k\|_{\mathsf Y}
\leq \frac{L_A}{2}\|h_k\|_{\mathsf X}^2
\leq \frac{L_AM^2}{2}\,
\tau^2(1+\tau)^{-2k-2}\|\rho_0\|_{\mathsf Y}^2.
\end{equation}

To control the accumulated error, set
$\theta_k:=T(\psi_k)-\bar p_k$.
The open-loop identity
$DT_{\psi_k}[h_k]=\bar p_{k+1}-\bar p_k$
implies $\theta_{k+1}=\theta_k+e_k$.
Since $\theta_0=0$, we have
$\theta_N=\sum_{k=0}^{N-1}e_k$.
Writing $C_{\mathrm{err}}:=L_AM^2/2$ and summing
\eqref{eq:error-geometric-one-step-bound} yields
\[
\begin{aligned}
\|\theta_N\|_{\mathsf Y}
&\leq C_{\mathrm{err}}\tau^2
\|\rho_0\|_{\mathsf Y}^2
\sum_{k=0}^{N-1}(1+\tau)^{-2k-2}\\
&=C_{\mathrm{err}}\frac{\tau}{2+\tau}
\bigl(1-(1+\tau)^{-2N}\bigr)
\|\rho_0\|_{\mathsf Y}^2\\
&\leq C_{\mathrm{err}}\frac{\tau}{2+\tau}
\|\rho_0\|_{\mathsf Y}^2.
\end{aligned}
\]
Taking the supremum over $N\geq0$ proves
\eqref{eq:holder-uniform-tracking}.

\noindent
Finally, the decomposition
$T(\psi_N)-\nu=\theta_N+(\bar p_N-\nu)$,
together with $\bar p_N-\nu=(1+\tau)^{-N}\rho_0$,
gives
\[
\|T(\psi_N)-\nu\|_{\mathsf Y}
\leq (1+\tau)^{-N}\|\rho_0\|_{\mathsf Y}
+C_{\mathrm{err}}\frac{\tau}{2+\tau}
\|\rho_0\|_{\mathsf Y}^2,
\]
which is \eqref{eq:holder-target-error}.
All constants are independent of $\tau$, $\rho_0$,
and the iteration index.
\end{proof}

\appendix

\bibliographystyle{unsrt}
\bibliography{ref}

\end{document}